\documentclass[onefignum,onetabnum]{siamart190516}

\usepackage{bm}
\usepackage{algorithm}
\usepackage{algpseudocode}
\usepackage{algorithmicx}
\usepackage{graphics}
\usepackage{subfigure}
\usepackage{amsfonts}
\usepackage{amssymb}
\usepackage{mathrsfs}
\allowdisplaybreaks
\usepackage{amsmath}  

\newtheorem{assumption}{Assumption}[section]

\newtheorem{remark}{Remark}[section]

\usepackage{overpic}
\def\ban{\begin{eqnarray*}}
	\def\ean{\end{eqnarray*}}
\def\bna{\begin{eqnarray}}
	\def\ena{\end{eqnarray}}
	
\usepackage{chngcntr} 
\counterwithout{algorithm}{section} 

\title{Distributed online learning with growing dimensions for stochastic large regression\thanks{This work was supported by the Natural Science Foundation
		of China under Grants No. 62603361 and No. T2293772. 
        \newline
			\hspace*{2em} 
				D. Gan is with the School of Mathematical Sciences, Fudan University, Shanghai 200433, China, and Laboratory of Mathematics for Nonlinear Science, Fudan University, Shanghai 200433, China, (e-mail: gandie@amss.ac.cn).
			Zhixin Liu is with State Key Laboratory of Mathematical Sciences, Academy of Mathematics and Systems Science, Chinese Academy of Sciences, Beijing 100190, China, (e-mail: lzx@amss.ac.cn).
	}}
	
	\author{Die Gan, Zhixin Liu\thanks{ Corresponding author.}}
	
\begin{document}

		\maketitle
		
		\begin{abstract}
			This paper investigates distributed online estimation and prediction for stochastic large regression models with infinitely many unknown parameters. Existing distributed estimation methods are mostly restricted to finite-dimensional parameter spaces, making them inadequate for emerging infinite-dimensional modeling applications.  To tackle this problem, this paper proposes a fully distributed online recursive least squares algorithm with dynamically growing model dimension. 
            The algorithm is constructed with an incrementally augmentable recursive structure, where new parameter dimensions  can be embedded in the recursion via zero-padded covariance and estimate augmentation without re-initialization. 
			The almost sure convergence of the algorithm is established under a weak cooperative non-persistent excitation condition, without requiring independence or stationarity assumptions on the regressors. 
			This condition integrates spatio-temporal network information and extends existing finite-dimensional results. The paper further  derives an almost sure upper bound on the accumulated prediction regret without any excitation conditions. It is proven that the average regret converges to zero almost surely for a broad class of dimension growth sequences, with the optimal finite-dimensional regret bound recovered as a special case. Finally, numerical simulations are conducted to verify the effectiveness of our theoretical findings.
		\end{abstract}
		
		




\section{Introduction}
Distributed estimation and adaptive learning for networked multi-agent systems have emerged as a critical research area  in modern control and signal processing, with wide-ranging applications including smart grid monitoring and multi-robot environment perception \cite{Simpson2016, Majcherczyk2021}. Unlike centralized frameworks that rely on a single fusion center, distributed approaches enable agents to cooperatively estimate unknown parameters using only local computations and information exchange with their neighbors. This distributed architecture offers  benefits in scalability, robustness to node and communication failures, and privacy preservation, making it essential for modern cyber-physical systems where centralized processing is impractical, bandwidth-limited, or privacy-prohibitive \cite{Li_2023}. 

A fundamental problem in this field is to identify the unknown parameters from historical trajectories  and 
delivering real-time prediction of a dynamical system's future states.  
In this work, we consider a set of  $n$ agents, where each agent $i$ observes noisy measurements of a signal pair $\{y_{k,i},\phi_{k,i}\}$ 
linked through a linear stochastic regression model of the form:
\begin{align}
	y^{\top}_{k+1,i} = \phi^{\top}_{k,i} \Theta + w^{\top}_{k+1,i}\label{system}
\end{align}
where $\Theta$ is a common unknown parameter and $\phi_{k,i}$ is the stochastic regression vector. 
The agents can communicate over a network, allowing information exchange with their neighbors. We focus on distributed online asymptotic estimation of  $\Theta$. 
Although, in principle, each agent could independently run a local algorithm to estimate  $\Theta$ from its own measurements, communication enables agents to leverage information 
collected by their neighbors. This cooperation leads to better convergence under excitation conditions that are strictly weaker than those required when relying solely on local data. 
Consequently, a distributed approach is particularly advantageous in scenarios where measurements are costly or difficult to obtain
in spatially distributed sensor systems \cite{Azzollini_2023}, with 
wireless environmental monitoring being a typical example.

Distributed RLS algorithms are widely studied for their fast convergence in estimating the unknown parameter \(\Theta\) in system (\ref{system}).
For instance, 
  \cite{yu2019} analyzed a diffusion RLS algorithm under the assumption that input regressors   
  $\phi_{k,i}$ are zero-mean and spatially independent.
\cite{Alireza2021} examined the mean-square convergence and stability of a diffusion RLS algorithm, assuming independent regression vectors 
$\phi_{k,i}$ 
with time-invariant covariance matrices.
\cite{Azzollini_2023} presented convergence analysis for a robust distributed RLS algorithm under a cooperative persistent excitation (PE) condition  which is equivalent to:  $\exists~ \alpha, T \geq 0, \quad \forall t \geq 0, \quad \sum_{k=t}^{t+T-1} \sum_{i =1}^n \phi_{k,i} \phi_{k,i}^\top \geq \alpha \bm I.$
In summary,  the majority of existing studies on  distributed RLS algorithm for stochastic dynamical systems require at least the persistence of excitation (PE) condition on the regression data, and most adopt restrictive independence or stationarity assumptions.  Although these idealized conditions simplify theoretical analysis, they are rarely satisfied or verifiable in many stochastic control systems with feedback signals  \cite{guo2020}. Recently,
some efforts have been devoted to relaxing the restrictive statistical assumptions  in distributed estimation for finite-dimensional systems \cite{Xie2021, Wang2021, Yan_2025, gan2026siam}.
Among these studies, \cite{Xie2021} derived the asymptotic convergence of distributed RLS estimators for finite-dimensional systems under the weakest possible cooperative non-PE condition on  signals. This condition guarantees consistent estimation even when individual agents are insufficiently excited, owing to the cooperation across the network.
The cooperative excitation concept has been further extended to various practical scenarios, including systems with unknown orders \cite{Gan202260}, colored noise environments \cite{xie2026}, and time-varying communication topologies \cite{Yan_2025, WANG2026}.

However, emerging challenges in machine learning,  system
identification, and adaptive control increasingly necessitate
handling high-dimensional, or even  infinite-dimensional unknown
parameters. Such models  arise in areas like the approximation
of infinite-width neural networks, nonparametric system
identification of distributed parameter systems, kernel-based
adaptive learning, and the modeling of emergent behavior in
large language models        \cite{Zhang_2025, Lv_2021, wei2022emergent, deep2025}. Traditional finite-dimensional distributed estimation methods including all the aforementioned works 
are inadequate in this infinite-dimensional context, as fixed-dimensional truncation introduces structural mismatch and non-vanishing steady-state errors.
For infinite-dimensional parameter identification, \cite{Sudipta2025} investigated the identification of transfer function coefficients in  single-input single-output infinite-dimensional systems, whereas \cite{Yin2022} proposed an infinite-dimensional sparse learning algorithm based on atomic norm regularization for linear system identification. 
\cite{zhang2022sieve} developed a sieve‑based stochastic gradient descent estimator  for online nonparametric regression over Sobolev ellipsoids using an increasing number of basis functions.
More closely aligned with the problem addressed in this paper are the works of \cite{DAI2025, DAI20251}, which first established the strong consistency of parameter estimates and the prediction performance of large models with saturated output observations.
These single-agent or centralized  results \cite{Sudipta2025,Yin2022,DAI2025, DAI20251} lay a theoretical foundation for the estimation of large regression models. However, they are not applicable to distributed networked systems, as they rely on a single fusion node with access to the global streaming dataset. Despite the growing interest in infinite-dimensional distributed optimization and control, the existing literature on distributed adaptive estimation for infinite-dimensional models is few.

Note that our previous works \cite{Gan2025CCC, wangying_IFAC_2026} have established the convergence of distributed RLS algorithms for infinite-dimensional stochastic large models.
However, they need to store a large amount of data  and recalculate all past data at each time step, which results in low computational efficiency.   Motivated by this limitation and inspired by \cite{DAI20251}, we  develop a fully distributed online  RLS algorithm  with growing model dimension to estimate the infinite unknown parameters 
for stochastic large regression models. The primary theoretical challenge lies in analyzing the convergence of recursive updates involving products of time-varying, non-independent random matrices with dynamically growing dimensions.  To tackle this challenge, we leverage  martingale theory, stochastic Lyapunov function methods, and algebraic graph theory for the theoretical analysis. The main contributions of this paper are summarized as follows:
\begin{itemize}
	\item  Compared with the  algorithms in \cite{Gan2025CCC, wangying_IFAC_2026}, we propose a fully online, and recursive distributed estimation algorithm  that automatically adjusts the parameter dimension as streaming data arrives.  The proposed algorithm consists of two main steps: 
    (1) a local RLS-based update step that manages dynamic dimension growth  using a truncation-augmentation mechanism, which initializes newly added parameter dimensions to zero without algorithm re-initialization; 
    and (2) a neighbor information fusion step based on the covariance intersection rule to ensure network consensus. This design  reduces computational complexity relative to prior growing batch-style methods (cf.,\cite{Gan2025CCC, wangying_IFAC_2026}),  making it suitable for long-term online learning.
	
	\item We establish almost sure convergence of the proposed algorithm without requiring the restrictive independence or stationarity assumptions on regressors in most existing literature. To this end, we introduce a cooperative non-PE condition tailored for infinite-dimensional large regression models. This condition extends the classical finite-dimensional cooperative non-PE result from \cite{Xie2021} and is much weaker than the standard PE condition (cf.,\cite{Azzollini_2023}). Also, it integrates both temporal and spatial network information to  capture multi-agent cooperative behavior, a feature absent in the single-agent conditions of \cite{DAI2025, DAI20251}.
	

	\item We derive an almost sure upper bound on the accumulated prediction regret without imposing any excitation conditions on the regressors. Moreover, we prove that the average regret converges to zero almost surely for a broad class of dimension growth sequences, including logarithmic and sublinear polynomial growth rates. These results provide rigorous theoretical guarantees for the proposed algorithm's online learning performance, and recover the classical optimal regret bound for finite-dimensional distributed RLS algorithms in \cite{Xie2021} as a special case, when the regressor dimension is bounded by a fixed constant.
\end{itemize}

The rest of the paper is organized as follows: Section \ref{sec:formulation} gives the problem formulation.		
Section \ref{sec:algorithm}  details the distributed online least squares algorithm with adaptive dimension growth. Section \ref{main_results} presents the main theoretical results on strong consistency and regret bounds. Numerical simulation results are provided in Section \ref{simulation}, with concluding remarks in Section \ref{section_conclusion}. 
	
\section{Problem Formulation}\label{sec:formulation}

\subsection{Matrix theory}
Let  $\mathbb{R}^m$ and $\mathbb{R}^{m\times n}$ denote the sets of $m$-dimensional real vectors and $m \times n$ real matrices, respectively. $\bm{I}_m$ and $\bm{0}_m$ represent the $m$-dimensional identity and zero matrices. For any matrix $A \in \mathbb{R}^{m \times n}$, its Euclidean norm is defined as $\|A\| \triangleq \sqrt{\lambda_{\max}(AA^\top)}$, where $(\cdot)^\top$ denotes matrix transposition, while $\lambda_{\max}(\cdot)$ and $\lambda_{\min}(\cdot)$ correspond to the largest and smallest eigenvalues, respectively. The matrix trace is denoted by $\operatorname{tr}(\cdot)$. The operators $\mathrm{col}(\cdot,\dots,\cdot)$ and $\mathrm{diag}(\cdot,\dots,\cdot)$ construct a column vector through vertical stacking and form a block-diagonal matrix, respectively. The Kronecker product is denoted by $\otimes$. For symmetric matrices $X, Y \in \mathbb{R}^{n \times n}$,  $X \geq Y$ and $X > Y$ indicate that $X-Y$ is positive semi-definite and positive definite, respectively (with analogous definitions for $\leq$ and $<$). Regarding matrix sequences $\{A_t\}$ and positive scalar sequences $\{a_t\}$, $A_t = O(a_t)$ means that $\|A_t\| \leq C a_t$ for all $t$ and some constant $C>0$, while $A_k = o(b_k)$ means $\lim_{k \to \infty} \|A_k\| / b_k = 0$.

\begin{lemma}[Matrix inversion formula \cite{Zielke1968}] \label{wl1}
	For matrices $ A$, $ B$,  $ C$ and  $ D$ of compatible dimensions, we have
	$( A+ B D C)^{-1}\!=\! A^{-1}- A^{-1} B( D^{-1}+ C A^{-1} B)^{-1} C  A^{-1},$ provided that the relevant matrices are invertible.
\end{lemma}		

\subsection{Communication topology}
Consider a network of $n$ agents modeled as a weighted directed graph $\mathcal{G}=(\mathcal{I},\mathcal{E},\mathcal{A})$, where $\mathcal{I}=\{1, 2,\dots, n\}$ represents the vertex set (agents), $\mathcal{E}\subseteq \mathcal{I}\times \mathcal{I}$ denotes the edge set, and $\mathcal{A}=[a_{ij}]_{n\times n}$ is the corresponding weighted adjacency matrix. An edge $(i, j)\in\mathcal{E}$, directed from agent $j$ to agent $i$, indicates that agent $i$ can obtain information from agent $j$. In this context, agent $j$ is referred to as an in-neighbor of agent $i$, and the set of all such neighbors for agent $i$ is denoted by $\mathcal{N}_{i}=\{j\in \mathcal{I} \mid (i,j)\in \mathcal{E}\}$. The elements of $\mathcal{A}$ are nonnegative, with $a_{ij}>0$ if and only if $(i, j)\in\mathcal{E}$.  Here we assume that $\mathcal{A}$ is a stochastic matrix.

A path of length $\ell$ in $\mathcal{G}$ consists of a sequence of vertices $\{v_{1},\dots,v_{\ell}\}$ such that $(v_{j},v_{j+1})\in\mathcal{E}$ for $1\leq j\leq\ell-1$. The digraph $\mathcal{G}$ is said to be strongly connected if there exists a directed path between every pair of distinct vertices. The diameter of $\mathcal{G}$, denoted by $D_{\mathcal{G}}$, is defined as the length of the maximum shortest path between any two vertices. Finally, $\mathcal{G}$ is called balanced if the in-degree and out-degree coincide for each vertex, i.e., $\sum_{j=1}^{n}a_{ji}=\sum_{j=1}^{n}a_{ij}=1$ for all $i=1,\dots,n$.

\subsection{System model}
This paper considers  a network comprising \(n\)  agents. At each time step \(k\), the observation of agent \(i\) is generated by the following stochastic large regression model:
\begin{align}
	y^{\top}_{k+1,i} = \phi^{\top}_{k,i} \Theta + w^{\top}_{k+1,i}, \quad k \geq 0,
	\label{model}
\end{align}
where \(y_{k,i} \in \mathbb{R}^m\) denotes the observation vector of agent \(i\) at time \(k\), \(\phi_{k,i} \in \mathbb{R}^{\infty}\) is an infinite-dimensional stochastic regressor with countably many components, \(w_{k+1,i} \in \mathbb{R}^m\) represents the measurement noise, and \(\Theta \in \mathbb{R}^{\infty \times m}\) is the unknown parameter matrix to be estimated. The matrix \(\Theta\) is assumed to satisfy the summability condition:
\begin{gather}
	\sum^{\infty}_{q=1} \| \Theta^{[q]} \| < \infty,
	\label{con}
\end{gather}
in which \(\Theta^{[q]} \in \mathbb{R}^m\) denotes the \(q\)-th row of \(\Theta\).

Note that model (\ref{model}) 
naturally includes the finite-dimensional case, as condition (\ref{con}) is automatically satisfied when $\Theta^{[q]}$
vanishes for sufficiently large  $q$. Moreover, the model encompasses a broad class of systems, including those with time-varying or expanding parameter dimensions.
It can be viewed as a single-hidden-layer neural network with infinite width, whose inner-layer parameters are pre-trained and fixed. This model also  provides a useful framework for transforming certain nonlinear or nonparametric problems into linear or parametric ones \cite{zhang2022sieve}.

We assume that the network communication topology of $n$ agents is  described by a weighted
digraph $\mathcal{G}=(\mathcal{I},\mathcal{E},\mathcal{A})$.
The goal of this paper is to design a distributed online estimation algorithm that enables all agents to cooperatively estimate the unknown infinite parameters $\Theta\in\mathbb{R}^{\infty\times m}$ using  the local information $\{y_{k,j},\phi_{k,j},k\geq 1\}_{j\in\mathcal{N}_{i}}$.

\section{Distributed estimation algorithm}\label{sec:algorithm}
This section presents the distributed algorithm design.
The key idea is to approximate the underlying infinite-dimensional parameter using a sequence of time-varying finite-dimensional truncations, whose effective dimension \(p_k\) grows gradually over time. Then, let $\{p_k\}$ be any non-decreasing sequence of positive integers and $p_k=O(k)$. 
Define the $p_k$- and $p_{k-1}$-dimensional regression vectors as
\begin{align}
	&\varphi_{k,i}=\text{col}\left\{\phi_{k,i}^{[1]},\phi_{k,i}^{[2]}\cdots, \phi_{k,i}^{[p_k]}\right\}\in\mathbb{R}^{p_k},\label{varphi}\\
	&\bar\varphi_{k,i}=\text{col}\left\{\phi_{k,i}^{[1]},\phi_{k,i}^{[2]}\cdots, \phi_{k,i}^{[p_{k-1}]}\right\}\in\mathbb{R}^{p_{k-1}},\label{phi}
\end{align}
where $\phi_{k,i}^{[q]}$ represents the $q$-th element of the vector $\phi_{k,i}$.
Built upon the above principle and inspired by  \cite{DAI20251}, we develop a distributed online least-squares algorithm with dynamically growing dimensions based on the regression vectors $\varphi_{k,i}$ and $\bar\varphi_{k,i}$. In this framework, the first \(p_{k-1}\) components of the parameter estimate are updated recursively from the previous step, while newly added dimensions are initialized to zero and  embedded into the existing recursive structure without full re-initialization. When \(p_k\) grows at a suitable speed (the specific choice of \(p_k\)  can be seen in Section \ref{main_results}), the algorithm has sufficient time to  learn the dominant parameter components in the leading dimensions through distributed iteration before expanding the model order.  The details can be seen in Algorithm \ref{algorithm}. The algorithm consists of two main steps. The first step performs a local recursive update extended from the standard RLS method, equipped with a dimension-adaptive augmentation scheme to handle dimension variations. This scheme automatically adjusts the update rules for the estimate and covariance matrix when the model dimension changes, avoiding a  restart of the estimation process. The second step fuses information from neighboring agents via a convex combination of local estimates and their inverse covariance matrices. This fusion strategy is based on the covariance intersection rule, which is a standard approach to achieving network-wide estimation consensus in distributed settings (cf.,\cite{Julier2001,Xie2021}). 

\begin{algorithm}
\caption{Distributed online LS algorithm with growing dimensions}\label{algorithm}
	Initialize: $P_{0,i}=I_{p_0}$, $\forall~ {\hat\Theta}_{0,i}\in\mathbb{R}^{p_0\times m}$.\\
	\textbf{ Step 1: Update auxiliary variable ${\bar\Theta}_{k,i}$ and ${\bar P}_{k,i}$ based  on $y_{k+1,i},$ ${\hat\Theta}_{k-1,i},$ $ P_{k-1,i},$ $\varphi_{k,i}$:}
	\begin{align}
		\bar{\Theta}_{k,i} 
		= &
		\mathbb{I}_{\{p_k = p_{k-1}\}}\left\{\!
		\hat{\Theta}_{k-1,i} \!+\! b_{k,i} P_{k-1,i} \varphi_{k,i} \left( y^{\top}_{k+1,i} \!-\! \bar\varphi_{k,i}^\top \hat{\Theta}_{k-1,i}\! \right)
		\right\} \nonumber\\
		+& \mathbb{I}_{\{p_k > p_{k-1}\}}\!
		\left\{\!
		\hat{\Theta}_{k-1,i}^* \!+\! b_{k,i} P_{k-1,i}^* \varphi_{k,i} \left( y^{\top}_{k+1,i} \!-\! \varphi_{k,i}^\top \hat{\Theta}^*_{k-1,i} \right)\!
		\right\}\!,
	\end{align}
	\begin{align}
		\bar{P}_{k,i} = &
		\mathbb{I}_{\{p_k = p_{k-1}\}}
		\left\{ 
		P_{k-1,i} - b_{k,i} P_{k-1,i} \varphi_{k,i} \varphi_{k,i}^\top P_{k-1,i}
		\right\} \nonumber\\
		& + \mathbb{I}_{\{p_k > p_{k-1}\}}
		\left\{ 
		P_{k-1,i}^* - b_{k,i} P_{k-1,i}^* \varphi_{k,i} \varphi_{k,i}^\top P_{k-1,i}^*
		\right\},\\
		b_{k,i} =& \frac{1}{1 + \varphi_{k,i}^{\top} \left[ \mathbb{I}_{\{p_k = p_{k-1}\}} P_{k-1,i} +  \mathbb{I}_{\{p_k > p_{k-1}\}} P_{k-1,i}^* \right] \varphi_{k,i}},
	\end{align}
	where
	\begin{align}
		&\hat{\Theta}_{k-1,i}^* \triangleq \mathrm{col}\{ \hat{\Theta}_{k-1,i}, \bm 0_{1\times m}, \cdots, \bm 0_{1\times m} \}\in \mathbb{R}^{p_k\times m},\\
		&P_{k-1,i}^* \triangleq
		\begin{pmatrix}
			P_{k-1,i} & 0 \\
			0 & \bm I_{p_k-p_{k-1}}
		\end{pmatrix}\in \mathbb{R}^{p_k\times p_k}.\label{9}
	\end{align}
	
	\textbf{ Step 2: Generate the matrix  $ P^{-1}_{k,i}$ and estimate ${\hat\Theta}_{k,i}$  by a convex combination of ${\bar\Theta}_{k,j}$ and ${\bar P}^{-1}_{k,j}$:}
	\begin{align*}
		P_{k,i}^{-1} = \sum_{j \in \mathcal{N}_{i}} a_{i j} \bar{P}_{k, j}^{-1} \in \mathbb{R}^{p_{k} \times p_{k}},~~
		\hat\Theta_{k, i} = P_{k i} \sum_{j \in \mathcal{N}_{i}} a_{i j} \bar{P}_{k ,j}^{-1} \bar{\Theta}_{k ,j} \in \mathbb{R}^{p_{k}\times m}.
	\end{align*}
\end{algorithm}


\begin{remark}
	Compared with our previous works in \cite{Gan2025CCC,wangying_IFAC_2026}, which require re-initialization and re-execution upon the arrival of new data, the proposed Algorithm \ref{algorithm} in this paper realizes a fully distributed online recursive estimation framework with significantly enhanced adaptability. Specifically, the algorithm can automatically adjust the dimension of the unknown parameter vector as new data arrive, enabling real-time parameter updates and output prediction without the need for re-initialization and re-execution. Moreover, Algorithm \ref{algorithm}  achieves an improvement in computational efficiency. The algorithms in \cite{Gan2025CCC,wangying_IFAC_2026} incur
	per-agent computational complexity  of $O(tp^3_t+t|\mathcal{N}_i|m p^2_t)$  at time $t$, while the proposed method in this paper reduces it to $O(p^3_t+|\mathcal{N}_i|m p^2_t)$ (see Table \ref{table1}), making it more suitable for long-term or large-scale adaptive learning scenarios.
\end{remark}

\begin{table}[!htb]
	\centering
	\caption{Computational complexity of several distributed algorithms for large regression models (per agent  at time $t$)}\label{table1}
	\begin{tabular}{lcc}
		\hline
		~& Algorithm type & Computational complexity   \\
		\hline
		Algorithm \ref{algorithm} (proposed)& Fully online & $O(p^3_t+|\mathcal{N}_i|m p^2_t)$ \\
		Algorithm in \cite{Gan2025CCC}& Not fully online& $O(tp^3_t+t|\mathcal{N}_i|m p^2_t)$  \\
		Algorithm in \cite{wangying_IFAC_2026} & Not fully online& $O(tp^3_t+t|\mathcal{N}_i|m p^2_t)$ \\
		\hline
	\end{tabular}\\
	\noindent \textit{Note:} $m$ is the dimension of observation vector $y_{k,i}$,
	$|\mathcal{N}_i|$ is the number of
	neighbors for agent $i$.
\end{table}

For convenience of analysis, we need to introduce the following 
notations:
\begin{align*}
	&{Y}_{k+1}^{\top}\triangleq{\rm col}\{{y}_{k+1,1}^{\top},\dots,{y}_{k+1,n}^{\top}\}, &(n\times m)\\
	&{\Phi}_{k}^{\top}\triangleq{\rm diag}\{{\varphi}_{k,1}^{\top},\dots,{\varphi}_{k,n}^{\top}\}, &(n\times np_k)\\
	&{\bar\Phi}_{k}^{\top}\triangleq{\rm diag}\{{\bar\varphi}_{k,1}^{\top},\dots,{\bar\varphi}_{k,n}^{\top}\}, &(n\times np_{k-1})\\
	&{W}_{k+1}^{\top}\triangleq{\rm col}\{{w}_{k+1,1}^{\top},\dots,{w}_{k+1,n}^{\top}\}, &(n\times m)\\
	&\Theta(k)\triangleq{\rm col}\{\Theta^{[1]},\cdots,\Theta^{[p_k]}\}, &(p_k\times m)\\
	&{\Gamma(k)}\triangleq{\rm col}\{\underbrace{{\Theta(k)},\dots,{\Theta(k)}}_{n}\}, &(np_k\times m)\\
	&{\hat\Theta}_{k}\triangleq{\rm col}\{{\hat\Theta}_{k,1},\dots,{\hat\Theta}_{k,n}\}, &(np_k\times m)\\
	&{\hat\Theta}^*_{k-1}\triangleq{\rm col}\{{\hat\Theta}^*_{k-1,1},\dots,{\hat\Theta}^*_{k-1,n}\}, &(np_k\times m)\\
	&\bar{{\Theta}}_{k}\triangleq{\rm col}\{\bar{{\Theta}}_{k,1},\dots,\bar{{\Theta}}_{k,n}\}, &(np_k\times m)\\
	&\widetilde{{\Theta}}_{k}\triangleq{\rm col}\{\widetilde{{\Theta}}_{k,1},\dots,\widetilde{{\Theta}}_{k,n}\}, &(np_k\times m)\\
	&~~~~~~~\text{where}~ \widetilde{{\Theta}}_{k,i}={\Theta(k)}-{\hat\Theta}_{k,i},\\
	&\widetilde{\bar{{\Theta}}}_{k}\triangleq{\rm col}\{\widetilde{\bar{{\Theta}}}_{k,1},\dots,\widetilde{\bar{{\Theta}}}_{k,n}\}, &(np_k\times m)\\
	&~~~~~~~\text{where}~ \widetilde{\bar{{\Theta}}}_{k,i}={\Theta(k)}-\bar{{\Theta}}_{k,i},\\
	&{b}_{k}\triangleq{\rm diag}\{b_{k,1},\dots,b_{k,n}\}, &(n\times n)\\
	&{c}_{k}\triangleq{b}_{k}\otimes I_{p_k}, &(np_k\times np_k)\\
	&{P}_{k}\triangleq{\rm diag}\{P_{k,1},\dots,P_{k,n}\}, &(np_k\times np_k)\\
	&\bar{{P}}_{k}\triangleq{\rm diag}\{\bar{P}_{k,1},\dots,\bar{P}_{k,n}\}, &(np_k\times np_k)\\
	&{\mathscr{A}_k}\triangleq\mathcal{A}\otimes I_{p_k}, &(np_k\times np_k)\\
	&\epsilon^{\top}_{k,i}=\phi^{\top}_{k,i}\Theta-\bar\varphi^{\top}_{k,i}\Theta(k-1)\!=\!\!\!\sum^{\infty}_{q=p_{k-1}+1}\!\!\!\!\phi^{[q]}_{k,i}\Theta^{[q]}, &(1\times m)\\
	&\epsilon^{\top}_k=\text{col}\{\epsilon^{\top}_{k,1},\cdots,\epsilon^{\top}_{k,n}\}. &(n\times m)\\
	&\epsilon'^{\top}_{k,i}=\phi^{\top}_{k,i}\Theta-\varphi^{\top}_{k,i}\Theta(k)=\sum^{\infty}_{q=p_{k}+1}\phi^{[q]}_{k,i}\Theta^{[q]}, &(1\times m)\\
	&\epsilon'^{\top}_k=\text{col}\{\epsilon'^{\top}_{k,1},\cdots,\epsilon'^{\top}_{k,n}\}. &(n\times m)
\end{align*}
By the definition of $\epsilon^{\top}_{k,i}$ and (\ref{model}), we have 
\begin{align}
	{y}_{k+1,i}^{\top}=\bar\varphi^{\top}_{k,i}\Theta(k-1)+w^{\top}_{k+1,i}+\epsilon^{\top}_{k,i}
	=\varphi^{\top}_{k,i}\Theta(k)+w^{\top}_{k+1,i}+\epsilon'^{\top}_{k,i}.\label{y1}
\end{align}
Thus we obtain the following compact form:
\begin{align}
	{Y}_{k+1}^{\top}=\bar\Phi^{\top}_{k}{\Gamma(k-1)}+W^{\top}_{k+1}+\epsilon^{\top}_{k}=\Phi^{\top}_{k}{\Gamma(k)}+W^{\top}_{k+1}+\epsilon'^{\top}_{k}
	.\label{1}
\end{align}
Since $\varphi_{k,i}^\top \hat{\Theta}^*_{k-1,i}=\bar\varphi_{k,i}^\top \hat{\Theta}_{k-1,i}$, then by Step 1 in Algorithm \ref{algorithm}, we have
\begin{align}
		\bar{\Theta}_{k} = &
		\mathbb{I}_{\{p_k = p_{k-1}\}}\left\{
		\hat{\Theta}_{k-1} + c_k P_{k-1} \Phi_{k} \left( Y^{\top}_{k+1} - \bar\Phi_{k}^\top \hat{\Theta}_{k-1} \right)
		\right\} \nonumber\\
		+ &\mathbb{I}_{\{p_k > p_{k-1}\}}
		\left\{
		\hat{\Theta}_{k-1}^* + c_k P_{k-1}^* \Phi_{k} \left( Y^{\top}_{k+1} - \Phi_{k}^\top \hat{\Theta}^*_{k-1} \right)
		\right\},\label{12}\\
		\bar{P}_{k} = &
		\mathbb{I}_{\{p_k = p_{k-1}\}}
		\left\{ 
		P_{k-1} - c_k P_{k-1} \Phi_{k} \Phi_{k}^\top P_{k-1}
		\right\} \nonumber\\
		+& \mathbb{I}_{\{p_k > p_{k-1}\}}
		\left\{ 
		P_{k-1}^* - c_k P_{k-1}^* \Phi_{k} \Phi_{k}^\top P_{k-1}^*
		\right\},\label{13}\\
		b_k =& \left\{I_{n} \!+\! \Phi_{k}^{\top} \left( \mathbb{I}_{\{p_k = p_{k-1}\}} P_{k-1} \!+\!  \mathbb{I}_{\{p_k > p_{k-1}\}} P_{k-1}^* \right) \Phi_{k}\right\}^{-1}.\label{14}
	\end{align}
	Substituting (\ref{1}) into the above equation, we can obtain the following error equation of Algorithm \ref{algorithm}.
	\begin{align}
		&\widetilde{\bar{{\Theta}}}_{k}=\Gamma(k)-\bar{{\Theta}}_{k}\nonumber\\
		=&\mathbb{I}_{\{p_k = p_{k-1}\}}\left\{\Gamma(k\!-\!1)\!-\!
		\hat{\Theta}_{k-1} \!- \!c_k P_{k-1} \Phi_{k} \!\left(\! Y^{\top}_{k+1} \!-\! \bar\Phi_{k}^\top \hat{\Theta}_{k-1} \right)\!
		\right\} \nonumber\\
		&+ \mathbb{I}_{\{p_k > p_{k-1}\}}
		\left\{\Gamma(k)\!-\!
		\hat{\Theta}_{k-1}^* \!-\!c_k P_{k-1}^* \Phi_{k} \!\left(\! Y^{\top}_{k+1} \!-\! \Phi_{k}^\top \hat{\Theta}^*_{k-1} \right)\!
		\right\}\nonumber\\
		=&	\mathbb{I}_{\{p_k = p_{k-1}\}}\left\{
		\widetilde{\Theta}_{k-1} - c_k P_{k-1} \Phi_{k} \left(  \bar\Phi_{k}^\top \widetilde{\Theta}_{k-1} +W^{\top}_{k+1}+\epsilon^{\top}_{k}\right)
		\right\} \nonumber\\
		&+ \mathbb{I}_{\{p_k > p_{k-1}\}}
		\Big\{\Gamma(k)\!-\!
		\hat{\Theta}_{k-1}^* \!-\!c_k P_{k-1}^* \Phi_{k} \left(  \Phi_{k}^\top (\Gamma(k)\!-\!
		\hat{\Theta}_{k-1}^*) +W^{\top}_{k+1}+\epsilon'^{\top}_{k} \right)
		\Big\}\nonumber\\
		=&	\mathbb{I}_{\{p_k = p_{k-1}\}}\Big\{\Big(\bm I_{p_{k-1}}-c_k P_{k-1} \Phi_{k}\Phi_{k}^\top\Big)
		\widetilde{\Theta}_{k-1} - c_k P_{k-1} \Phi_{k} \left(   W^{\top}_{k+1}+\epsilon^{\top}_{k}\right)
		\Big\} \nonumber\\
		&+ \mathbb{I}_{\{p_k > p_{k-1}\}}\Big\{\Big(\bm I_{p_{k}}-c_k P^*_{k-1} \Phi_{k}\Phi_{k}^\top\Big)(\Gamma(k)\!-\!
		\hat{\Theta}_{k-1}^*)
	- c_k P^*_{k-1} \Phi_{k} \left(   W^{\top}_{k+1}+\epsilon'^{\top}_{k}\right)
		\Big\} \nonumber\\
		=&	\mathbb{I}_{\{p_k = p_{k-1}\}}\Big\{\Big(P^{-1}_{k-1}+\Phi_k\Phi^{\top}_k\Big)^{-1}P^{-1}_{k-1}
		\widetilde{\Theta}_{k-1} - c_k P_{k-1} \Phi_{k} \left(   W^{\top}_{k+1}+\epsilon^{\top}_{k}\right)
		\Big\} \nonumber\\
		&+ \mathbb{I}_{\{p_k > p_{k-1}\}}\Big\{\Big(P^{*-1}_{k-1}+\Phi_k\Phi^{\top}_k\Big)^{-1}P^{*-1}_{k-1}(\Gamma(k)\!-\!
		\hat{\Theta}_{k-1}^*)
	 - c_k P^*_{k-1} \Phi_{k} \left(   W^{\top}_{k+1}+\epsilon'^{\top}_{k}\right)
		\Big\}, \label{bartheta}
	\end{align}
	where Lemma \ref{wl1} is used in last equation.
	
	Note that by (\ref{13})-(\ref{14}) and Lemma \ref{wl1}, we have 
	\begin{align}
		{\bar P}^{-1}_{k}=\mathbb{I}_{\{p_k = p_{k-1}\}}P^{-1}_{k-1}+\mathbb{I}_{\{p_k > p_{k-1}\}}P^{*-1}_{k-1}+\Phi_k\Phi^{\top}_k.\label{barpni}
	\end{align}

	By Step 2 in Algorithm \ref{algorithm}, for each agent $i\in\{1,2,...,n\}$
	$$
	P_{k,i}\sum_{j\in \mathcal{N}_{i}}a_{ij}{\bar P}^{-1}_{k,j}{\Theta}(k)={\Theta}(k).
	$$
	Thus we obtain that 
	\begin{align}
		P_{k}\mathscr{A}_k{\bar P}^{-1}_{k}\Gamma(k)
		=\text{col}\Big\{P_{k,1}\sum_{j\in \mathcal{N}_{1}}a_{1j}{\bar P}^{-1}_{k,j}{\Theta}(k),\cdots,P_{k,n}\sum_{j\in \mathcal{N}_{n}}a^{(k)}_{nj}{\bar P}^{-1}_{k,j}{\Theta}(k)\Big\}=\Gamma(k).
	\end{align}
	Hence by Step 2 in Algorithm \ref{algorithm} and (\ref{bartheta}), we have
	\begin{align}
		\widetilde{\Theta}_{k}
		=&\Gamma(k)-\hat\Theta_{k}
		=
		\Gamma(k)-P_{k}{\mathscr{A}_k}{\bar P}^{-1}_{k}\bar{{\Theta}}_{k}
		=P_{k}{\mathscr{A}_k}{\bar P}^{-1}_{k}\widetilde{\bar{{\Theta}}}_{k}\nonumber\\
		=&\mathbb{I}_{\{p_k = p_{k-1}\}}\Big\{P_{k}{\mathscr{A}_k}P^{-1}_{k-1}
		\widetilde{\Theta}_{k-1} - P_{k}\mathscr{A}_k{\bar P}^{-1}_{k}c_k P_{k-1} \Phi_{k} \left(   W^{\top}_{k+1}+\epsilon^{\top}_{k}\right)
		\Big\} \nonumber\\
		&+ \mathbb{I}_{\{p_k > p_{k-1}\}}\Big\{P_{k}\mathscr{A}_kP^{*-1}_{k-1}(\Gamma(k)\!-\!
		\hat{\Theta}_{k-1}^*)
		\nonumber\\
		&\hskip 2cm - P_{k}\mathscr{A}_k{\bar P}^{-1}_{k}c_k P^*_{k-1} \Phi_{k} \left(   W^{\top}_{k+1}+\epsilon'^{\top}_{k}\right)
		\Big\} \nonumber\\
		\triangleq&\mathbb{I}_{\{p_k = p_{k-1}\}}M_k+\mathbb{I}_{\{p_k > p_{k-1}\}}N_k.
		\label{error_equation}
	\end{align}

	\section{The Main results}\label{main_results}
	In this section, we will establish the strong consistency of Algorithm \ref{algorithm} under a cooperative excitation condition of the regressors, and will analyze  the prediction error without relying on any excitation conditions.  To begin with, we introduce some basic assumptions about noise and graph.
	\subsection{Assumptions}
	\begin{assumption}\label{AssuSPR}
		Suppose that for any $i=1,\dots,n$, each element of $\varphi_{k,i}$  is  $\mathscr{F}_k$-measurable, and the noise sequence $\{{w}_{k,i},\mathcal{F}_{k}\}$ is a martingale difference sequence, and there exists $\beta>2$ such that
		$$
		\sup_{k\geq 0}\mathbb{E}[\|{w}_{k+1,i}\|^{\beta}\vert\mathcal{F}_{k}]<\infty,~\text{a.s.,}$$
		where $\{\mathcal{F}_{k}\}$ is a sequence of nondecreasing $\sigma$-algebras.
	\end{assumption}
	{\begin{remark}
			Assumption \ref{AssuSPR}	is the standard condition for the noise in the theoretical analysis. It is clear that the i.i.d. zero-mean bounded or Gaussian noise independent of the regression signals can
			satisfy
			Assumption \ref{AssuSPR}.
		\end{remark}

		\begin{assumption}\label{AssuGrap}
			$\mathcal{G}$ is strongly connected and balanced.
		\end{assumption}

		\begin{remark}\label{RemaGrap}
			By Assumption \ref{AssuGrap}, we know that when $\ell$ is no less than the diameter of the graph $\mathcal{G}$, i.e., $D_{\mathcal{G}}$, each element of the $\ell$-th power of the matrix $\mathcal{A}$ shall be positive. Note that switching  graphs that satisfy certain connectivity assumptions can also achieve the convergence of the algorithm. The essential idea is the same except that one exploits the property of joint connectivity instead.
		\end{remark}

		\subsection{Key lemmas}
		
		For convenience, we denote the following notations:
		\begin{align}\label{nota}
			\begin{aligned}
				&\bar{\bar{\Theta}}(k)=\text{col}\{\Theta^{[p_{k-1}+1]},\cdots,\Theta^{[p_{k}]}\}\in\mathbb{R}^{(p_k-p_{k-1})\times m},\\
				&\bar{\bar{\Theta}}'_{k}\triangleq\text{col}\{\underbrace{\bar{\bar{\Theta}}(k),\dots,\bar{\bar{\Theta}}(k)}_{n}\}\in\mathbb{R}^{[n(p_k-p_{k-1})]\times m},\\
				&\bar{\bar{\varphi}}_{k,i}=\text{col}\{\phi^{[p_{k-1}+1]}_{k,i},\cdots,\phi^{[p_{k}]}_{k,i}\}\in\mathbb{R}^{p_k-p_{k-1}},\\
				&\bar{\bar{\Phi}}_k=\text{diag}\{\bar{\bar{\varphi}}_{k,1},\cdots,\bar{\bar{\varphi}}_{k,n}\}\in\mathbb{R}^{n(p_k-p_{k-1})\times n}.
			\end{aligned}
		\end{align}

			We first present a key property regarding the convex combination of positive semi-definite matrices.
		\begin{lemma}  \cite{Xie2021} \label{LemmCoro41}
			For any adjacency matrix $\mathcal{A}=\{a_{ij}\}\in\mathbb{R}^{n\times n}$ of a weight-balanced digraph, denote $\mathscr{A}_k=\mathcal{A}\otimes I_{p_k}$. Then for any $k\geq 1$, we have
			\begin{align*}
				\mathscr{A}_k^{\top}\bar{{P}}_{k}^{-1}\mathscr{A}_k\leq {P}_{k}^{-1},
			\end{align*}
			and
			\begin{align*}
				\mathscr{A}_k^{\top}{P}_{k}\mathscr{A}_k\leq \bar{{P}}_{k}.
			\end{align*}
		\end{lemma}		
		
		By Lemma \ref{LemmCoro41}, we have the following result for the stochastic Lyapunov function $\operatorname{tr}(	\widetilde{\Theta}^{\top}_{t}{P}_{t}^{-1}\widetilde{\Theta}_{t})$.
		\begin{lemma}\label{lemma4.2}
			Consider model (\ref{model})-(\ref{con})  and Algorithm \ref{algorithm}, 
			we have  the following property for the estimation error:
			\begin{align}
				\operatorname{tr}(	\widetilde{\Theta}^{\top}_{t}{P}_{t}^{-1}\widetilde{\Theta}_{t})
				\leq 
				\operatorname{tr}(	\widetilde{\Theta}^{\top}_{0}{P}_{0}^{-1}\widetilde{\Theta}_{0})+\sum^{t}_{k=1}R_k+\sum^{t}_{k=1}\mathbb{I}_{\{p_k > p_{k-1}\}}S_k, ~\text{a.s.,}
			\end{align}
			where
			\begin{align}	
				R_k		\triangleq 
				&-\operatorname{tr}\left(\widetilde{\Theta}^{\top}_{k-1}\bar\Phi_{k}b_k\bar\Phi^{\top}_{k}\widetilde{{\Theta}}_{k-1}
				\right)
				-2\operatorname{tr}\Big(\widetilde{\Theta}^{\top}_{k-1} \bar\Phi_{k}b_k[W^{\top}_{k+1}+\epsilon'^{\top}_{k}]
				\Big)\nonumber\\
				&+\operatorname{tr}\Big([W_{k+1}+\epsilon'_{k}]b_k\bar\Phi^{\top}_{k}{P}_{k-1}\bar\Phi_{k}[W^{\top}_{k+1}+\epsilon'^{\top}_{k}]
				\Big)
			\end{align}
			and 		
			\begin{align}	
				S_k	\triangleq& 
				\operatorname{tr}(\bar{\bar{\Theta}}'^{\top}_{k}\bar{\bar{\Theta}}'_{k})
				-2\operatorname{tr}(\widetilde{\Theta}^{\top}_{k-1}\bar\Phi_{k}b_k\bar{\bar{\Phi}}^{\top}_k\bar{\bar{\Theta}}'_{k})
				-\operatorname{tr}(\bar{\bar{\Theta}}'^{\top}_{k}\bar{\bar{\Phi}}_kb_k\bar{\bar{\Phi}}^{\top}_k\bar{\bar{\Theta}}'_{k})\nonumber\\
				&-2\operatorname{tr}(\bar{\bar{\Theta}}'^{\top}_{k}\bar{\bar{\Phi}}_kb_k[W^{\top}_{k+1}+\epsilon'^{\top}_{k}])+\operatorname{tr}([W_{k+1}+\epsilon'_{k}]b_k{\bar{\bar\Phi}}_{k}^{\top}{\bar{\bar\Phi}}_{k}[W^{\top}_{k+1}+\epsilon'^{\top}_{k}]).
			\end{align}		
		\end{lemma}	
		\begin{proof}
			By (\ref{error_equation}), we have
			\begin{align}
				V_{k}=&\widetilde{\Theta}^{\top}_{k}{P}_{k}^{-1}\widetilde{\Theta}_{k}\nonumber\\
				=&(\mathbb{I}_{\{p_k = p_{k-1}\}}M_k+\mathbb{I}_{\{p_k > p_{k-1}\}}N_k)^{\top}{P}_{k}^{-1}\cdot(\mathbb{I}_{\{p_k = p_{k-1}\}}M_k+I_{\{p_k> p_{k-1}\}}N_k)\nonumber\\
				=&\mathbb{I}_{\{p_k = p_{k-1}\}}M^{\top}_k{P}_{k}^{-1}M_k+\mathbb{I}_{\{p_k > p_{k-1}\}}N^{\top}_k{P}_{k}^{-1}N_k.
				\label{v1}
			\end{align}	
			In the following, we estimate $\mathbb{I}_{\{p_k = p_{k-1}\}}M^{\top}_k{P}_{k}^{-1}M_k$ and $\mathbb{I}_{\{p_k > p_{k-1}\}}N^{\top}_k{P}_{k}^{-1}N_k$, respectively.
			
			\textbf{(1) Compute $\mathbb{I}_{\{p_k = p_{k-1}\}}\operatorname{tr}(M^{\top}_k{P}_{k}^{-1}M_k$).}
			
			By the definition of $M_k$ in (\ref{error_equation}) and Lemma \ref{LemmCoro41}, we have 
			\begin{align}
				&\mathbb{I}_{\{p_k = p_{k-1}\}}M^{\top}_k{P}_{k}^{-1}M_k\nonumber\\
				=&\mathbb{I}_{\{p_k = p_{k-1}\}}\Big\{{
					\widetilde{\Theta}^{\top}_{k-1}P^{-1}_{k-1}\mathscr{A}^{\top}_k}P_{k}-\left(W_{k+1}+\epsilon_{k}\right) \Phi^{\top}_{k}P_{k-1}c_k {\bar P}^{-1}_{k}\mathscr{A}^{\top}_kP_{k} 
				\Big\} \nonumber\\
				&\cdot\Big\{{\mathscr{A}_k}P^{-1}_{k-1}
				\widetilde{\Theta}_{k-1}- \mathscr{A}_k{\bar P}^{-1}_{k}c_k P_{k-1} \Phi_{k} \left(   W^{\top}_{k+1}+\epsilon^{\top}_{k}\right)
				\Big\}\nonumber\\
				=&\mathbb{I}_{\{p_k = p_{k-1}\}}\Big\{
				\widetilde{\Theta}^{\top}_{k-1}P^{-1}_{k-1}-\left(W_{k+1}+\epsilon_{k}\right) \Phi^{\top}_{k}P_{k-1}c_k {\bar P}^{-1}_{k}
				\Big\} \nonumber\\
				&\cdot\mathscr{A}^{\top}_kP_{k}{\mathscr{A}_k}\Big\{P^{-1}_{k-1}
				\widetilde{\Theta}_{k-1}-{\bar P}^{-1}_{k}c_k P_{k-1} \Phi_{k} \left(   W^{\top}_{k+1}+\epsilon^{\top}_{k}\right)
				\Big\}\nonumber\\
				\leq&\mathbb{I}_{\{p_k = p_{k-1}\}}\Big\{
				\widetilde{\Theta}^{\top}_{k-1}P^{-1}_{k-1}-\left(W_{k+1}+\epsilon_{k}\right) \Phi^{\top}_{k}P_{k-1}c_k {\bar P}^{-1}_{k}
				\Big\} \nonumber\\
				&\cdot\bar P_{k}\Big\{P^{-1}_{k-1}
				\widetilde{\Theta}_{k-1}-{\bar P}^{-1}_{k}c_k P_{k-1} \Phi_{k} \left(   W^{\top}_{k+1}+\epsilon^{\top}_{k}\right)
				\Big\}.\nonumber
			\end{align}
			Thus it can be obtained by taking the trace operation on both sides of the above equation, 
			\begin{align}	
				&\mathbb{I}_{\{p_k = p_{k-1}\}}\operatorname{tr}(M^{\top}_k{P}_{k}^{-1}M_k)\nonumber\\
				\leq&\mathbb{I}_{\{p_k = p_{k-1}\}}\Bigg\{	\operatorname{tr}\left(\widetilde{\Theta}^{\top}_{k-1}{P}_{k-1}^{-1}{\bar P}_{k}{P}^{-1}_{k-1}\widetilde{{\Theta}}_{k-1}\right)\nonumber\\
				&-2\operatorname{tr}\Big(\widetilde{\Theta}^{\top}_{k-1}{P}_{k-1}^{-1} c_k{P}_{k-1}\Phi_{k}[W^{\top}_{k+1}+\epsilon^{\top}_{k}]
				\Big)\nonumber\\
				&+\operatorname{tr}\Big([W_{k+1}+\epsilon_{k}]\Phi^{\top}_{k}{P}_{k-1}c_k{\bar P}^{-1}_{k}c_k{P}_{k-1}\Phi_{k}[W^{\top}_{k+1}+\epsilon^{\top}_{k}]
				\Big)\Bigg\}.\label{v2}
			\end{align}
			Now, we proceed to compute the right-hand side (RHS) of
			(\ref{v2}) term by term. Firstly, we know that from (\ref{13})
			\begin{align}	
				&\mathbb{I}_{\{p_k = p_{k-1}\}}\operatorname{tr}\left(\widetilde{\Theta}^{\top}_{k-1}{P}_{k-1}^{-1}\bar{P}_{k}{P}^{-1}_{k-1}\widetilde{{\Theta}}_{k-1}\right)\nonumber\\
				=&\mathbb{I}_{\{p_k = p_{k-1}\}}\operatorname{tr}\Big(\widetilde{\Theta}^{\top}_{k-1}{P}_{k-1}^{-1}\Big({P}_{k-1}-c_k{P}_{k-1}\Phi_{k}\Phi^{\top}_{k}P_{k-1}\Big){P}^{-1}_{k-1}\widetilde{{\Theta}}_{k-1}\Big)\nonumber\\
				=&\mathbb{I}_{\{p_k = p_{k-1}\}}\operatorname{tr}\left(\widetilde{\Theta}^{\top}_{k-1}{P}_{k-1}^{-1}\widetilde{\Theta}_{k-1}\right)-\mathbb{I}_{\{p_k = p_{k-1}\}}\operatorname{tr}\left(\widetilde{\Theta}^{\top}_{k-1}\Phi_{k}b_k\Phi^{\top}_{k}\widetilde{{\Theta}}_{k-1}
				\right),\label{v3}
			\end{align}
			where the last inequality holds since we have
			\begin{align}\label{27}
				c_kP_{k-1}=P_{k-1}c_k,~\Phi^{\top}_{k}c_k=b_k\Phi^{\top}_{k}, ~c_k\Phi_{k}=\Phi_{k}b_k.
			\end{align}
			Then for the second term on the RHS of (\ref{v2}), we have  by (\ref{27})
			\begin{align}
				&\mathbb{I}_{\{p_k = p_{k-1}\}}\operatorname{tr}\Big(\widetilde{\Theta}^{\top}_{k-1}{P}_{k-1}^{-1} c_k{P}_{k-1}\Phi_{k}[W^{\top}_{k+1}+\epsilon^{\top}_{k}]
				\Big)\nonumber\\
				=&\mathbb{I}_{\{p_k = p_{k-1}\}}\operatorname{tr}\Big(\widetilde{\Theta}^{\top}_{k-1} \Phi_{k}b_k[W^{\top}_{k+1}+\epsilon^{\top}_{k}]
				\Big).\label{v4}   
			\end{align}
			For the last term on the RHS of (\ref{v2}), by Lemma \ref{LemmCoro41}, (\ref{14}) and  (\ref{barpni}), we compute it as follows:
			\begin{align}
				&\mathbb{I}_{\{p_k = p_{k-1}\}}\operatorname{tr}\Big([W_{k+1}+\epsilon_{k}]\Phi^{\top}_{k}{P}_{k-1}c_k{\bar P}^{-1}_{k}c_k{P}_{k-1}\Phi_{k}[W^{\top}_{k+1}+\epsilon^{\top}_{k}]
				\Big)\nonumber\\
				=&\mathbb{I}_{\{p_k = p_{k-1}\}}\operatorname{tr}\Big([W_{k+1}+\epsilon_{k}]\Phi^{\top}_{k}{P}_{k-1}c_k[P^{-1}_{k-1}+\Phi_{k}\Phi^{\top}_{k}]c_k{P}_{k-1}\Phi_{k}[W^{\top}_{k+1}+\epsilon^{\top}_{k}]
				\Big)
				\nonumber\\
				=&\mathbb{I}_{\{p_k = p_{k-1}\}}\operatorname{tr}\Big([W_{k+1}+\epsilon_{k}]b^2_k\Phi^{\top}_{k}{P}_{k-1}\Phi_{k}[W^{\top}_{k+1}+\epsilon^{\top}_{k}]
				\Big)	\nonumber\\
				&+\mathbb{I}_{\{p_k = p_{k-1}\}}\operatorname{tr}\Big([W_{k+1}+\epsilon_{k}]b^2_k\Phi^{\top}_{k}{P}_{k-1}\Phi_{k}\Phi^{\top}_{k} {P}_{k-1}\Phi_{k}[W^{\top}_{k+1}+\epsilon^{\top}_{k}]
				\Big)\nonumber\\	
				=&\mathbb{I}_{\{p_k = p_{k-1}\}}\operatorname{tr}\Big([W_{k+1}+\epsilon_{k}]b_k\Phi^{\top}_{k}{P}_{k-1}\Phi_{k}[W^{\top}_{k+1}+\epsilon^{\top}_{k}]
				\Big).\label{v5}
			\end{align}
			Thus by (\ref{v2}), (\ref{v3}), (\ref{v4}) and (\ref{v5}), we can obtain that
			\begin{align}
				&\mathbb{I}_{\{p_k = p_{k-1}\}}\operatorname{tr}(M^{\top}_k{P}_{k}^{-1}M_k)\nonumber\\
				\leq& \mathbb{I}_{\{p_k = p_{k-1}\}}\operatorname{tr}\left(\widetilde{\Theta}^{\top}_{k-1}{P}_{k-1}^{-1}\widetilde{\Theta}_{k-1}\right)-\mathbb{I}_{\{p_k = p_{k-1}\}}\operatorname{tr}\left(\widetilde{\Theta}^{\top}_{k-1}\Phi_{k}b_k\Phi^{\top}_{k}\widetilde{{\Theta}}_{k-1}
				\right)\nonumber\\
				&-2\mathbb{I}_{\{p_k = p_{k-1}\}}\operatorname{tr}\Big(\widetilde{\Theta}^{\top}_{k-1} \Phi_{k}b_k[W^{\top}_{k+1}+\epsilon^{\top}_{k}]
				\Big)\nonumber\\
				&+\mathbb{I}_{\{p_k = p_{k-1}\}}\operatorname{tr}\Big([W_{k+1}+\epsilon_{k}]b_k\Phi^{\top}_{k}{P}_{k-1}\Phi_{k}[W^{\top}_{k+1}+\epsilon^{\top}_{k}]
				\Big).\label{M}
			\end{align}
			
			\textbf{(2) Compute $\mathbb{I}_{\{p_k > p_{k-1}\}}\operatorname{tr}(N^{\top}_k{P}_{k}^{-1}N_k$).}
			
			Note that the definition of $N_k$ is similar to $M_k$ (see (\ref{error_equation})). Thus following the proof of equation (\ref{M}), we can obtain that
			\begin{align}
				&\mathbb{I}_{\{p_k > p_{k-1}\}}\operatorname{tr}(N^{\top}_k{P}_{k}^{-1}N_k)\nonumber\\
				\leq& \mathbb{I}_{\{p_k > p_{k-1}\}}\operatorname{tr}\left({(\Gamma(k)\!-\!
					\hat{\Theta}_{k-1}^*)^{\top}}{P}_{k-1}^{*-1}(\Gamma(k)\!-\!
				\hat{\Theta}_{k-1}^*)\right)\nonumber\\
				&-\mathbb{I}_{\{p_k > p_{k-1}\}}\operatorname{tr}\left({(\Gamma(k)\!-\!
					\hat{\Theta}_{k-1}^*)^{\top}}\Phi_{k}b_k\Phi^{\top}_{k}(\Gamma(k)\!-\!
				\hat{\Theta}_{k-1}^*)
				\right)\nonumber\\
				&-2\mathbb{I}_{\{p_k > p_{k-1}\}}\operatorname{tr}\Big({(\Gamma(k)\!-\!
					\hat{\Theta}_{k-1}^*)^{\top}} \Phi_{k}b_k[W^{\top}_{k+1}+\epsilon'^{\top}_{k}]
				\Big)\nonumber\\
				&+\mathbb{I}_{\{p_k > p_{k-1}\}}\operatorname{tr}\Big([W_{k+1}+\epsilon'_{k}]b_k\Phi^{\top}_{k}{P}^*_{k-1}\Phi_{k}[W^{\top}_{k+1}+\epsilon'^{\top}_{k}]
				\Big).\label{N}
			\end{align}
			Then by the definition of $\Gamma(k)$, $\hat{\Theta}_{k-1}^*$ and (\ref{nota}), we have for $p_k>p_{k-1}$
			\begin{align}
				\Gamma(k)\!-\!
				\hat{\Theta}_{k-1}^*&=\text{col}\{\Theta(k)-\hat{\Theta}_{k-1,1}^*,...,\Theta(k)-\hat{\Theta}_{k-1,n}^*\}\nonumber\\
				&=\text{col}\left\{\begin{pmatrix}  \widetilde{{\Theta}}_{k-1,1}\\ \bar{\bar{\Theta}}_{k}\end{pmatrix},...,\begin{pmatrix}  \widetilde{{\Theta}}_{k-1,n}\\ \bar{\bar{\Theta}}_{k}\end{pmatrix}\right\}.\label{err}
			\end{align}
			By the definition of ${P}^*_{k-1}$, (\ref{9}) and (\ref{err}), it is clear that
			\begin{align}
				&\mathbb{I}_{\{p_k > p_{k-1}\}}\operatorname{tr}\left({(\Gamma(k)\!-\!
					\hat{\Theta}_{k-1}^*)^{\top}}{P}_{k-1}^{*-1}(\Gamma(k)\!-\!
				\hat{\Theta}_{k-1}^*)\right)\nonumber\\
				=&\mathbb{I}_{\{p_k > p_{k-1}\}}\left\{\operatorname{tr}\left(\widetilde{\Theta}^{\top}_{k-1}{P}_{k-1}^{-1}\widetilde{\Theta}_{k-1}\right)+\operatorname{tr}(\bar{\bar{\Theta}}'^{\top}_{k}\bar{\bar{\Theta}}'_{k})\right\}.\label{29}
			\end{align} 
			Also, by the definition of $\Phi_k$, (\ref{nota}) and (\ref{err}), we have	
			\begin{align}
				&\mathbb{I}_{\{p_k > p_{k-1}\}}(\Gamma(k)\!-\!
				\hat{\Theta}_{k-1}^*)^{\top}\Phi_{k}\nonumber\\
				=&\mathbb{I}_{\{p_k > p_{k-1}\}}\left(\text{col}\left\{\begin{pmatrix}  \widetilde{{\Theta}}_{k-1,1}\\ \bar{\bar{\Theta}}_{k}\end{pmatrix},...,\begin{pmatrix}  \widetilde{{\Theta}}_{k-1,n}\\ \bar{\bar{\Theta}}_{k}\end{pmatrix}\right\}\right)^{\top}\text{diag}\left\{\begin{pmatrix}  \bar\varphi_{k,1}\\ {\bar{\bar\varphi}}_{k,1}\end{pmatrix},...,\begin{pmatrix}  \bar\varphi_{k,n}\\ {\bar{\bar\varphi}}_{k,n}\end{pmatrix}\right\}
				\nonumber\\
				=&\mathbb{I}_{\{p_k > p_{k-1}\}}\text{col}\Big\{\widetilde{{\Theta}}^{\top}_{k-1,1}\bar\varphi_{k,1}+\bar{\bar{\Theta}}^{\top}_{k}{\bar{\bar\varphi}}_{k,1},\cdots,\widetilde{{\Theta}}^{\top}_{k-1,n}\bar\varphi_{k,n}+\bar{\bar{\Theta}}^{\top}_{k}{\bar{\bar\varphi}}_{k,n}\Big\}\nonumber\\
				=&\mathbb{I}_{\{p_k > p_{k-1}\}}(\widetilde{{\Theta}}^{\top}_{k-1}\bar\Phi_{k}+\bar{\bar{\Theta}}'^{\top}_{k}\bar{\bar{\Phi}}_k).\label{30}
			\end{align}		
			By the definitions of ${P}^*_{k-1}$	and  $\Phi_{k}$, it follows that from  (\ref{nota})
			\begin{align}		
				&\mathbb{I}_{\{p_k > p_{k-1}\}}\Phi^{\top}_{k}{P}^*_{k-1}\Phi_{k}\nonumber\\
				=&\mathbb{I}_{\{p_k > p_{k-1}\}}\text{diag}\left\{\begin{pmatrix}  \bar\varphi_{k,1}\\ {\bar{\bar\varphi}}_{k,1}\end{pmatrix}^{\top},...,\begin{pmatrix}  \bar\varphi_{k,n}\\ {\bar{\bar\varphi}}_{k,n}\end{pmatrix}^{\top}\right\}\nonumber\\
				&\cdot\text{diag}\left\{\begin{pmatrix}
					P_{k-1,1} & 0 \\
					0 & \bm I_{p_k-p_{k-1}}
				\end{pmatrix},\cdots,\begin{pmatrix}
					P_{k-1,n} & 0 \\
					0 & \bm I_{p_k-p_{k-1}}
				\end{pmatrix}\right\}\nonumber\\
				&\cdot\text{diag}\left\{\begin{pmatrix}  \bar\varphi_{k,1}\\ {\bar{\bar\varphi}}_{k,1}\end{pmatrix},...,\begin{pmatrix}  \bar\varphi_{k,n}\\ {\bar{\bar\varphi}}_{k,n}\end{pmatrix}\right\}
				\nonumber\\
				=&\mathbb{I}_{\{p_k > p_{k-1}\}}\text{diag}\Big\{\bar\varphi_{k,1}^{\top}P_{k-1,1}\bar\varphi_{k,1}+{\bar{\bar\varphi}}_{k,1}^{\top}{\bar{\bar\varphi}}_{k,1}, \cdots,\bar\varphi_{k,n}^{\top}P_{k-1,n}\bar\varphi_{k,n}+{\bar{\bar\varphi}}_{k,n}^{\top}{\bar{\bar\varphi}}_{k,n}
				\Big\}\nonumber\\
				=&\mathbb{I}_{\{p_k > p_{k-1}\}}(\bar\Phi_{k}^{\top}P_{k-1}\bar\Phi_{k}+{\bar{\bar\Phi}}_{k}^{\top}{\bar{\bar\Phi}}_{k}
				).\label{31}
			\end{align}		
			Substituting (\ref{29})-(\ref{31}) into (\ref{N}),  we can obtain that
			\begin{align}
				&\mathbb{I}_{\{p_k > p_{k-1}\}}\operatorname{tr}(N^{\top}_k{P}_{k}^{-1}N_k)\nonumber\\
				\leq&\mathbb{I}_{\{p_k > p_{k-1}\}}\left\{\operatorname{tr}\left(\widetilde{\Theta}^{\top}_{k-1}{P}_{k-1}^{-1}\widetilde{\Theta}_{k-1}\right)+\operatorname{tr}(\bar{\bar{\Theta}}'^{\top}_{k}\bar{\bar{\Theta}}'_{k})\right\}\nonumber\\
				-&\mathbb{I}_{\{p_k > p_{k-1}\}}\operatorname{tr}\left((\widetilde{{\Theta}}^{\top}_{k-1}\bar\Phi_{k}+\bar{\bar{\Theta}}'^{\top}_{k}\bar{\bar{\Phi}}_k)b_k(\bar\Phi^{\top}_{k}\widetilde{{\Theta}}_{k-1}+\bar{\bar{\Phi}}^{\top}_k\bar{\bar{\Theta}}'_{k})
				\right)\nonumber\\
				&-2\mathbb{I}_{\{p_k > p_{k-1}\}}\operatorname{tr}\Big((\widetilde{{\Theta}}^{\top}_{k-1}\bar\Phi_{k}+\bar{\bar{\Theta}}'^{\top}_{k}\bar{\bar{\Phi}}_k)b_k[W^{\top}_{k+1}+\epsilon'^{\top}_{k}]
				\Big)\nonumber\\
				&+\mathbb{I}_{\{p_k > p_{k-1}\}}\operatorname{tr}\Big([W_{k+1}+\epsilon'_{k}] b_k(\bar\Phi_{k}^{\top}P_{k-1}\bar\Phi_{k}+{\bar{\bar\Phi}}_{k}^{\top}{\bar{\bar\Phi}}_{k}
				)[W^{\top}_{k+1}+\epsilon'^{\top}_{k}]
				\Big)\nonumber\\
				=&\mathbb{I}_{\{p_k > p_{k-1}\}}\operatorname{tr}\left(\widetilde{\Theta}^{\top}_{k-1}{P}_{k-1}^{-1}\widetilde{\Theta}_{k-1}\right)+\mathbb{I}_{\{p_k > p_{k-1}\}}(R_k+S_k),\label{32}
			\end{align}
			Note that  $\epsilon^{\top}_{k}=\epsilon'^{\top}_{k}$ and $\Phi_k=\bar \Phi_k$ hold when $p_k=p_{k-1}$. Then from (\ref{M}), (\ref{32}) and (\ref{v1}), we can obtain that 
			\begin{align}
				\operatorname{tr}(V_{k})\!=&\mathbb{I}_{\{p_k = p_{k-1}\}}\!\operatorname{tr}(M^{\top}_k{P}_{k}^{-1}M_k)\!+\!\mathbb{I}_{\{p_k > p_{k-1}\}}\operatorname{tr}(N^{\top}_k{P}_{k}^{-1}N_k)\nonumber\\
				\leq&\operatorname{tr}(	V_{k-1})\!+\!\mathbb{I}_{\{p_k > p_{k-1}\}}R_k+\mathbb{I}_{\{p_k > p_{k-1}\}}(R_k+S_k)\nonumber\\
				=&\operatorname{tr}(V_{k-1})+R_k+\mathbb{I}_{\{p_k > p_{k-1}\}}S_k.
			\end{align}					
			Summing from $k=1$ to $t$ yields the desired result, which completes the proof of this theorem.
		\end{proof}

		\begin{remark}
			Lemma \ref{lemma4.2} establishes an almost-sure upper bound on $\operatorname{tr}(\!	\widetilde{\Theta}^{\top}_{t}\!{P}_{t}^{-1}\!\widetilde{\Theta}_{t}\!)$, which concerns the estimation error under growing model dimensions. 	Intuitively, the total error up to time $t$ is bounded by the sum of three components: the initial estimation error with dimension $p_0$, a recursive learning term $R_k$ with dimension $p_{k-1}$, and an expansion cost term $S_k$
			that is activated only when the parameter dimension increases, i.e., $p_k>p_{k-1}$. Then, if the dimension expansions occur sufficiently slowly, the error reduction achieved through $R_k$
			may compensate for the cost $S_k$
			incurred during each expansion.
		\end{remark}

		Thus based on Lemma \ref{lemma4.2}, we can further obtain the following result.
		\begin{lemma}\label{lemma4.3}
			Consider model (\ref{model})-(\ref{con})  and Algorithm \ref{algorithm}. Under Assumption \ref{AssuSPR},
			we have the estimation property:
			\begin{align}
				&\operatorname{tr}(	\widetilde{\Theta}^{\top}_{t}{P}_{t}^{-1}\widetilde{\Theta}_{t})+\left(\frac{1}{4}-o(1)\right)\sum^{t}_{k=1}\operatorname{tr}\left(\widetilde{\Theta}^{\top}_{k-1}\bar\Phi_{k}b_k\bar\Phi^{\top}_{k}\widetilde{{\Theta}}_{k-1}
				\right)
				\nonumber\\
				=&O\left(p_t\log^+\{ \lambda_{\max}(P^{-1}_t)\}+p_t+\sum^n_{i=1}\sum^{t}_{k=1}\|\epsilon_{k,i}\|^2\right),
			\end{align}
				where 
			$$r_t=1+\sum^n_{i=1}\sum^t_{k=1}\|\widetilde\varphi_{k,i}(t)\|^2$$
			with $\widetilde\varphi_{k,i}(t)\!\triangleq\! \text{col}\left\{\phi_{k,i}^{[1]},\phi_{k,i}^{[2]}\cdots, \phi_{k,i}^{[p_k]},\underbrace{0,\cdots,0}_{p_t-p_k}\right\}\in\mathbb{R}^{p_t},$
			and $\log^+\{x\}=\max\{\log x,0\}$.
		\end{lemma}

		\begin{proof}
			Based on Lemma \ref{lemma4.2}, we  should estimate  the matrices  $\sum^{t}_{k=1}R_k$  and $\sum^{t}_{k=1}I_\{p_k>p_{k-1}\}S_k$. 
			%
			By  the properties of the trace and the martingale estimation theorem (Theorem 2.8 in \cite{Guo1991}) and $b_{k,i}\leq 1$, we can  obtain the following estimation for any $\delta>0$,
			\begin{align}
				&\sum^{t}_{k=1}\operatorname{tr}\Big(\widetilde{\Theta}^{\top}_{k-1}
				\bar\Phi_{k}b_kW^{\top}_{k+1}\Big)
				=\sum^n_{i=1}\sum^{t}_{k=1}\operatorname{tr}\Big(\widetilde{\Theta}^{\top}_{k-1,i}
				\bar\varphi_{k,i}b_{k,i}w^{\top}_{k+1,i}\Big)\nonumber\\
				=&\sum^n_{i=1}\sum^{t}_{k=1}b_{k,i}\bar\varphi^{\top}_{k,i}\widetilde{\Theta}_{k-1,i}w_{k+1,i}
				=O\left(\sum^n_{i=1}\left\{\sum^{t}_{k=1}b_{k,i}\|\bar\varphi^{\top}_{k,i}\widetilde{\Theta}_{k-1,i}\|^2\right\}^{\frac{1}{2}+\delta}\right)\nonumber\\
				=&O\left(\sum^n_{i=1}\left\{o\left(\sum^{t}_{k=1}b_{k,i}\|\bar\varphi^{\top}_{k,i}\widetilde{\Theta}_{k-1,i}\|^2\right)+O(1)\right\}\right)\nonumber\\
				=&o\left(\sum^{t}_{k=1}\operatorname{tr}\Big(\widetilde{\Theta}^{\top}_{k-1}
				\bar\Phi_{k}b_k	\bar\Phi^{\top}_{k}\widetilde{\Theta}_{k-1}
				\Big)\right)+O(1).
				\label{47}
			\end{align}
			
			By the fact $\lambda_{\max}(b_k)\leq 1$ , it is clear that 
			\begin{align}
				-2\sum^{t}_{k=1}\operatorname{tr}\Big(\widetilde{\Theta}^{\top}_{k-1}
				\bar\Phi_{k}b_k\epsilon'^{\top}_{k}\Big)
				\leq& \frac{1}{2}\sum^{t}_{k=1}\operatorname{tr}\left(\widetilde{\Theta}^{\top}_{k-1}\bar\Phi_{k}b_k\bar\Phi^{\top}_{k}\widetilde{{\Theta}}_{k-1}
				\right)\!+2\!\sum^{t}_{k=1}\operatorname{tr}\left(\epsilon'^{\top}_{k}\epsilon'_{k}\right)\nonumber\\
				=& \frac{1}{2}\sum^{t}_{k=1}\operatorname{tr}\left(\widetilde{\Theta}^{\top}_{k-1}\bar\Phi_{k}b_k\bar\Phi^{\top}_{k}\widetilde{{\Theta}}_{k-1}
				\right)\!\!+2\sum^n_{i=1}\!\sum^{t}_{k=1}\|\epsilon'_{k,i}\|^2.\label{new5}
			\end{align}
			By (\ref{barpni}) and following the proof  of Lemma 4.4 in \cite{Xie2021} (see (33) in \cite{Xie2021}), we have for $p_k=p_{k-1}$,
			\begin{align}
				\lambda_{\max}(b_k\bar\Phi^{\top}_{k}{P}_{k-1}\bar\Phi_{k})
				\leq\frac{{\rm det}\left({P}_{k}^{-1}\right)-{\rm det}\left({P}_{k-1}^{-1}\right)}{{\rm det}\left({P}_{k}^{-1}\right)}.\label{new7}
			\end{align}
			and for $p_k>p_{k-1}$,
			\begin{align}
				\lambda_{\max}(b_k\Phi^{\top}_{k}{P}^*_{k-1}\Phi_{k})
				\leq\frac{{\rm det}\left({P}_{k}^{-1}\right)-{\rm det}\left({P}_{k-1}^{*-1}\right)}{{\rm det}\left({P}_{k}^{-1}\right)}
				=\frac{{\rm det}\left({P}_{k}^{-1}\right)-{\rm det}\left({P}_{k-1}^{-1}\right)}{{\rm det}\left({P}_{k}^{-1}\right)}.\label{new8}
			\end{align}
			Thus by (\ref{31}), (\ref{new7}) and (\ref{new8}), we can obtain that
			\begin{align}
				&\mathbb{I}_{\{p_k = p_{k-1}\}}	\lambda_{\max}(b_k\bar\Phi^{\top}_{k}{P}_{k-1}\bar\Phi_{k})
				+\mathbb{I}_{\{p_k > p_{k-1}\}}\lambda_{\max}(b_k\Phi^{\top}_{k}{P}^*_{k-1}\Phi_{k})\nonumber\\
				=&\mathbb{I}_{\{p_k = p_{k-1}\}}	\lambda_{\max}(b_k\bar\Phi^{\top}_{k}{P}_{k-1}\bar\Phi_{k})
				+\mathbb{I}_{\{p_k > p_{k-1}\}}\lambda_{\max}(b_k\bar\Phi^{\top}_{k}{P}_{k-1}\bar\Phi_{k}+b_k{\bar{\bar\Phi}}_{k}^{\top}{\bar{\bar\Phi}}_{k})\nonumber\\
				\leq&\frac{{\rm det}\left({P}_{k}^{-1}\right)-{\rm det}\left({P}_{k-1}^{-1}\right)}{{\rm det}\left({P}_{k}^{-1}\right)}\leq 1.\label{new9}
			\end{align}
			Then it follows from (\ref{new9}) that
			\begin{align}
				&\sum^{t}_{k=1}\operatorname{tr}\Big([W_{k+1}\!+\!\epsilon'_{k}] 
				b_k\bar\Phi^{\top}_{k}{P}_{k-1}\bar\Phi_{k}[W^{\top}_{k+1}\!+\!\epsilon'^{\top}_{k}]
				\Big)\nonumber\\
				\leq& 2\sum^{t}_{k=1}\operatorname{tr}\Big(W_{k+1}
				b_k\bar\Phi^{\top}_{k}{P}_{k-1}\bar\Phi_{k}W^{\top}_{k+1}
				\Big)+ 2\sum^{t}_{k=1}\operatorname{tr}\Big(\epsilon'_{k}
				b_k\bar\Phi^{\top}_{k}{P}_{k-1}\bar\Phi_{k}\epsilon'^{\top}_{k}
				\Big)\nonumber\\
				\leq& 2\sum^{t}_{k=1}\operatorname{tr}\Big(W_{k+1}
				b_k\bar\Phi^{\top}_{k}{P}_{k-1}\bar\Phi_{k}W^{\top}_{k+1}
				\Big)+ 2\sum^{t}_{k=1}\lambda_{\max}(b_k\bar\Phi^{\top}_{k}{P}_{k-1}\bar\Phi_{k})\operatorname{tr}\Big(\epsilon'_{k}
				\epsilon'^{\top}_{k}
				\Big)\nonumber\\
				\leq& 2\sum^{t}_{k=1}\operatorname{tr}\Big(W_{k+1}
				b_k\bar\Phi^{\top}_{k}{P}_{k-1}\bar\Phi_{k}W^{\top}_{k+1}
				\Big)+ 2\sum^n_{i=1}\!\sum^{t}_{k=1}\|\epsilon'_{k,i}\|^2.\label{36}
			\end{align}
			Similar to 	(\ref{47}) and (\ref{new5}) and (\ref{36}), for the  terms  in $\sum^{t}_{k=1}\mathbb{I}_{\{p_k > p_{k-1}\}}S_k$ defined in Lemma \ref{lemma4.2}, we have the following estimates:
			\begin{align}
				\sum^{t}_{k=1}\mathbb{I}_{\{p_k > p_{k-1}\}}\operatorname{tr}\Big(\bar{\bar{\Theta}}'^{\top}_{k}\bar{\bar{\Phi}}_kb_kW^{\top}_{k+1}\Big)
				=o\left(\sum^{t}_{k=1}\mathbb{I}_{\{p_k > p_{k-1}\}}\operatorname{tr}\Big(\bar{\bar{\Theta}}'^{\top}_{k}\bar{\bar{\Phi}}_kb_k\bar{\bar{\Phi}}^{\top}_k\bar{\bar{\Theta}}'_{k}
				\Big)\right)+O(1),\label{52}
			\end{align}
			and
			\begin{align}
				&-2\sum^{t}_{k=1}\mathbb{I}_{\{p_k > p_{k-1}\}}\operatorname{tr}\Big(\bar{\bar{\Theta}}'^{\top}_{k}\bar{\bar{\Phi}}_kb_k\epsilon'^{\top}_{k}\Big)\nonumber\\
				\leq& \frac{1}{2}\sum^{t}_{k=1}\mathbb{I}_{\{p_k > p_{k-1}\}}\operatorname{tr}\Big(\bar{\bar{\Theta}}'^{\top}_{k}\bar{\bar{\Phi}}_kb_k\bar{\bar{\Phi}}^{\top}_k\bar{\bar{\Theta}}'_{k}
				\Big)\!\!+2\sum^n_{i=1}\!\sum^{t}_{k=1}\mathbb{I}_{\{p_k > p_{k-1}\}}\|\epsilon'_{k,i}\|^2,
			\end{align}
			and 	
			\begin{align}
				&\sum^{t}_{k=1}\mathbb{I}_{\{p_k > p_{k-1}\}}\operatorname{tr}\Big([W_{k+1}\!+\!\epsilon'_{k}] 
				b_k{\bar{\bar\Phi}}_{k}^{\top}{\bar{\bar\Phi}}_{k}[W^{\top}_{k+1}\!+\!\epsilon'^{\top}_{k}]
				\Big)\nonumber\\
				\leq& 2\sum^{t}_{k=1}\mathbb{I}_{\{p_k > p_{k-1}\}}\operatorname{tr}\Big(W_{k+1}
				b_k{\bar{\bar\Phi}}_{k}^{\top}{\bar{\bar\Phi}}_{k}W^{\top}_{k+1}
				\Big)+ 2\sum^n_{i=1}\!\sum^{t}_{k=1}\mathbb{I}_{\{p_k > p_{k-1}\}}\|\epsilon'_{k,i}\|^2.\label{54}
			\end{align}
			By  the properties of the trace, we have
			\begin{align}
				&-2\sum^{t}_{k=1}\mathbb{I}_{\{p_k > p_{k-1}\}}\operatorname{tr}(\widetilde{\Theta}^{\top}_{k-1}\bar\Phi_{k}b_k\bar{\bar{\Phi}}^{\top}_k\bar{\bar{\Theta}}'_{k})\nonumber\\
				\leq & \frac{1}{4}\sum^{t}_{k=1}\operatorname{tr}\left(\widetilde{\Theta}^{\top}_{k-1}\bar\Phi_{k}b_k\bar\Phi^{\top}_{k}\widetilde{{\Theta}}_{k-1}\right)+4\sum^{t}_{k=1}\mathbb{I}_{\{p_k > p_{k-1}\}}\operatorname{tr}(\bar{\bar{\Theta}}'^{\top}_{k}\bar{\bar{\Phi}}_kb_k\bar{\bar{\Phi}}^{\top}_k\bar{\bar{\Theta}}'_{k}).\label{56}
			\end{align}
			Thus by Lemma \ref{lemma4.2}, (\ref{47}), (\ref{new5}), (\ref{36}), (\ref{52})-(\ref{56}), we have 
			\begin{align}
				&\operatorname{tr}(	\widetilde{\Theta}^{\top}_{t}{P}_{t}^{-1}\widetilde{\Theta}_{t})+\left(\frac{1}{4}-o(1)\right)\sum^{t}_{k=1}\operatorname{tr}\left(\widetilde{\Theta}^{\top}_{k-1}\bar\Phi_{k}b_k\bar\Phi^{\top}_{k}\widetilde{{\Theta}}_{k-1}
				\right)
				\nonumber\\
				\leq &2\sum^{t}_{k=1}\operatorname{tr}\Big(W_{k+1}
				b_k\bar\Phi^{\top}_{k}{P}_{k-1}\bar\Phi_{k}W^{\top}_{k+1}
				\Big)+\!
				2\sum^{t}_{k=1}\mathbb{I}_{\{p_k > p_{k-1}\}}\operatorname{tr}\Big(W_{k+1}
				b_k{\bar{\bar\Phi}}_{k}^{\top}{\bar{\bar\Phi}}_{k}W^{\top}_{k+1}	\Big)	\nonumber\\
				&+\sum^{t}_{k=1}\mathbb{I}_{\{p_k > p_{k-1}\}}	\operatorname{tr}(\bar{\bar{\Theta}}'^{\top}_{k}\bar{\bar{\Theta}}'_{k})+O\left(\sum^{t}_{k=1}\mathbb{I}_{\{p_k > p_{k-1}\}}\operatorname{tr}(\bar{\bar{\Theta}}'^{\top}_{k}\bar{\bar{\Phi}}_kb_k\bar{\bar{\Phi}}^{\top}_k\bar{\bar{\Theta}}'_{k})\right)\nonumber\\
				&+8\sum^n_{i=1}\!\sum^{t}_{k=1}\|\epsilon'_{k,i}\|^2\!+\!O(1).\label{48}
			\end{align}
			In the following, we estimate the right-hand side (RHS) of (\ref{48}) term by term.
			By (\ref{new9}) and $P^{-1}_{0}=I_{p_0}$, we have
			\begin{align}
				&\sum^{t}_{k=1}\lambda_{\max}(b_k\bar\Phi^{\top}_{k}{P}_{k-1}\bar\Phi_{k})
				\leq \sum^{t}_{k=1}\frac{{\rm det}\left({P}_{k}^{-1}\right)-{\rm det}\left({P}_{k-1}^{-1}\right)}{{\rm det}\left({P}_{k}^{-1}\right)}\nonumber\\
				=&\sum^{t}_{k=1}\int^{{\rm det}\left({P}_{k}^{-1}\right)}_{{\rm det}\left({P}_{k-1}^{-1}\right)}\frac{1}{{\rm det}\left({P}_{k}^{-1}\right)}dx
				\leq \int^{{\rm det}\left({P}_{t}^{-1}\right)}_{{\rm det}\left({P}_{0}^{-1}\right)}\frac{1}{x}dx\nonumber\\
				=&\log({\rm det}\left({P}_{t}^{-1}\right))-\log({\rm det}\left({P}_{0}^{-1}\right))
				=O\left(p_t\log^+\{ \lambda_{\max}(P^{-1}_t)\}\right).\label{41}
			\end{align}
			Similar to (\ref{41}), we have
			\begin{align}
				\sum^{t}_{k=1}\sum^{t}_{k=1}\mathbb{I}_{\{p_k > p_{k-1}\}}\lambda_{\max}(b_k{\bar{\bar\Phi}}_{k}^{\top}{\bar{\bar\Phi}}_{k})
				=O\left(p_t\log^+\{ \lambda_{\max}(P^{-1}_t)\}\right).\label{45}
			\end{align}
			
			Thus from martingale estimation theorem, Assumption \ref{AssuSPR} and (\ref{41}) it follows that
			\begin{align}
				&	\sum^{t}_{k=1}\operatorname{tr}\Big(W_{k+1}
				b_k\bar\Phi^{\top}_{k}{P}_{k-1}\bar\Phi_{k}W^{\top}_{k+1}
				\Big)\nonumber\\	
				\leq &\sum^{t}_{k=1}\!\!\left(\lambda_{\max}(b_k\bar\Phi^{\top}_{k}{P}_{k-1}\bar\Phi_{k})\!\sum^n_{i=1}\!\|w_{k+1,i}\|^2\right)\nonumber\\	
				\leq &\Bigg\{\sum^{t}_{k=1}\Big(\lambda_{\max}(b_k\bar\Phi^{\top}_{k}{P}_{k-1}\bar\Phi_{k})\sum^n_{i=1}\!(\|w_{k+1,i}\|^2-\mathbb{E}[\|w_{k+1,i}\|^{2}|\mathscr{F}_k]) \Big)\nonumber\\	
				+&\sum^{t}_{k=1}\lambda_{\max}(b_k\bar\Phi^{\top}_{k}{P}_{k-1}\bar\Phi_{k})\sum^n_{i=1}\mathbb{E}[\|w_{k+1,i}\|^{2}|\mathscr{F}_k]) 
				\Bigg\}\nonumber\\	
				=&[O(1)+o(1)]\left(\sum^{t}_{k=1}\lambda_{\max}(b_k\bar\Phi^{\top}_{k}{P}_{k-1}\bar\Phi_{k})\right)+O(1)\nonumber\\	
				=&O\left(p_t\log^+\{ \lambda_{\max}(P^{-1}_t)\}\right)
				\!+\!O(1).\label{42}
			\end{align}
			Similar to (\ref{42}), we can obtain that from (\ref{45})
			\begin{align}
				&2\sum^{t}_{k=1}\mathbb{I}_{\{p_k > p_{k-1}\}}\operatorname{tr}\Big(W_{k+1}
				b_k{\bar{\bar\Phi}}_{k}^{\top}{\bar{\bar\Phi}}_{k}W^{\top}_{k+1}	\Big)\nonumber\\
				=&[O(1)+o(1)]\left(\sum^{t}_{k=1}\sum^{t}_{k=1}\mathbb{I}_{\{p_k > p_{k-1}\}}\lambda_{\max}(b_k{\bar{\bar\Phi}}_{k}^{\top}{\bar{\bar\Phi}}_{k})\right)+O(1)\nonumber\\
				=&O\left(p_t\log^+\{ \lambda_{\max}(P^{-1}_t)\}\right)
				\!+\!O(1).
			\end{align}
			By (\ref{con}) and the definition of $\bar{\bar{\Theta}}'$ in (\ref{nota}), we know that 
			\begin{align}
				\sum^{t}_{k=1}\mathbb{I}_{\{p_k > p_{k-1}\}}	\operatorname{tr}(\bar{\bar{\Theta}}'^{\top}_{k}\bar{\bar{\Theta}}'_{k})=O(p_t).
			\end{align}
			By the definitions of $\bar{\bar{\Theta}}'$, $\bar{\bar{\Phi}}$, $\epsilon'_{k,i}$,$\epsilon_{k,i}$ and $b_k\leq \bm I_n$, we have
			\begin{align}
				&\left(\sum^{t}_{k=1}\mathbb{I}_{\{p_k > p_{k-1}\}}\operatorname{tr}(\bar{\bar{\Theta}}'^{\top}_{k}\bar{\bar{\Phi}}_kb_k\bar{\bar{\Phi}}^{\top}_k\bar{\bar{\Theta}}'_{k})\right)
				+\sum^n_{i=1}\!\sum^{t}_{k=1}\|\epsilon'_{k,i}\|^2\nonumber\\
				\leq&\sum^{n}_{i=1}\sum^{t}_{k=1}\sum^{p_k}_{q=p_{k-1}+1}\|\phi^{[q]}_{k,i}\Theta^{[q]}\|^2+\sum^{n}_{i=1}\sum^{t}_{k=1}\sum^{\infty}_{q=p_{k}+1}\|\phi^{[q]}_{k,i}\Theta^{[q]}\|^2\nonumber\\
				=&\sum^{n}_{i=1}\sum^{t}_{k=1}\sum^{\infty}_{q=p_{k-1}+1}\|\phi^{[q]}_{k,i}\Theta^{[q]}\|^2=\sum^n_{i=1}\!\sum^{t}_{k=1}\|\epsilon_{k,i}\|^2.\label{64}
			\end{align}
			Thus substituting (\ref{42})-(\ref{64}) into (\ref{48}), we can get the desired result.
		\end{proof}
		
		\subsection{Parameter convergence}
		 Based on Lemma \ref{lemma4.3}, we provide a theoretical upper bound for the parameter estimation error induced by Algorithm \ref{algorithm}.
		\begin{theorem}\label{th1}
			Consider model (\ref{model})-(\ref{con})  and Algorithm \ref{algorithm}. Under Assumptions \ref{AssuSPR} and \ref{AssuGrap},
			we have the following  property for the estimation error:
			\begin{align}
				\|{\hat\Theta}'_{t,i}-\Theta\|^2
				=O\!\left(\frac{1}{\lambda_{\min}(t)}\left[p_t\log r_t\!+\sum^n_{i=1}\sum^{t}_{k=1}\|\epsilon_{k,i}\|^2\right]+\gamma_t\right) \text{a.s.,}\label{th}
			\end{align}
			where
			\begin{align*}
				&{\hat\Theta}'_{t,i}\triangleq\text{col}\{{\hat\Theta}_{t,i},\bm 0_m,\bm 0_m,...\}\in\mathbb{R}^{\infty\times m},\\
				&r_t=1+\sum^n_{i=1}\sum^t_{k=1}\|\widetilde\varphi_{k,i}(t)\|^2,~		\gamma_t\triangleq\left(\sum^{\infty}_{q=p_t+1}\left\|\Theta^{[q]}\right\|\right)^2,\\
				&\widetilde\varphi_{k,i}(t)\triangleq \text{col}\left\{\phi_{k,i}^{[1]},\phi_{k,i}^{[2]}\cdots, \phi_{k,i}^{[p_k]},\underbrace{0,\cdots,0}_{p_t-p_k}\right\}\in\mathbb{R}^{p_t},\\
				&\lambda_{\min}(t)\!=\!\lambda_{\min}\left\{I_{p_t}+\sum^n_{i=1}\sum ^{t-D_{\mathcal{G}}}_{k=1} \widetilde\varphi_{k,i}(t){ \widetilde\varphi}^{\top}_{k,i}(t)\!\right\}
			\end{align*}
			satisfying $\lambda_{\min}(t)\xrightarrow[{t\rightarrow\infty}] {} \infty$.	
		\end{theorem}

		\begin{proof}	
			In the following, we estimate the matrix $P^{-1}_{t}$.
			First, we  prove the following equation  holds by induction:
			\begin{align}
				\text{vec}\{P^{-1}_{t}\}\!=\!\text{vec}\{\bm I_{np_t}\}\!+\!\sum^t_{k=1}\mathscr{A}_t^{t-k}\text{vec}\{\widetilde\Phi_k(t)\widetilde\Phi^{\top}_k(t)\},\label{vec}
			\end{align}
			where $\text{vec}\{\cdot\}$ denotes the operator that stacks the blocks of a block diagonal matrix on top of each other, and $$\widetilde\Phi_k(t)\triangleq\text{diag}\{\widetilde\varphi_{k,1}(t),\cdots,\widetilde\varphi_{k,n}(t)\}\in\mathbb{R}^{np_t\times n}.$$
			If 	equation (\ref{vec}) holds for $t-1$, then we have
			\begin{align}
				\text{vec}\{P^{-1}_{t-1}\}
				\!=\!\text{vec}\{\bm I_{np_{t-1}}\}
				+\sum^{t}_{k=1}\mathscr{A}_{t-1}^{t-1-k}\text{vec}\{\widetilde\Phi_k(t-1)\widetilde\Phi^{\top}_k(t-1)\}.\label{vec1}
			\end{align}
			This is equivalent to the following equation holding for every  $i\in\{1,...,n\}$:
			\begin{align}
				P^{-1}_{t-1,i}=\bm I_{p_{t-1}}+\sum^{t-1}_{k=1}\sum^n_{j=1}a^{(t-1-k)}_{ij}\widetilde\varphi_{k,j}(t-1)\widetilde\varphi^{\top}_{k,j}(t-1),\label{vec2}
			\end{align}
			where $a^{(t-1-k)}_{ij}$ is the $i$-th row, $j$-column element of the matrix $\mathcal{A}^{t-1-k}$ with $\mathcal{A}$ being the adjacency matrix.
			By Step 2 in Algorithm \ref{algorithm}, we know that 
			\begin{align}
				\text{vec}\{P^{-1}_{t}\}={\mathscr{A}_t}\text{vec}\{\bar P^{-1}_{t}\}.
			\end{align}
			Thus combining this with (\ref{barpni}) and (\ref{vec1}), we have
			\begin{align}
				\text{vec}\{P^{-1}_{t}\}
				=&{\mathscr{A}_t}\text{vec}\{\mathbb{I}_{\{p_t = p_{t-1}\}}P^{-1}_{t-1}+\mathbb{I}_{\{p_t > p_{t-1}\}}P^{*-1}_{t-1}+\Phi_t\Phi^{\top}_t\}\nonumber\\
				=&{\mathscr{A}_t}\Big\{\mathbb{I}_{\{p_t = p_{t-1}\}}\text{vec}\{P^{-1}_{t-1}\}+\mathbb{I}_{\{p_t > p_{t-1}\}}\text{vec}\{P^{*-1}_{t-1}\}+\text{vec}\{\Phi_t\Phi^{\top}_t\}\Big\}\nonumber\\
				=&{\mathscr{A}_t}\Big\{\mathbb{I}_{\{p_t = p_{t-1}\}}\Big(\text{vec}\{\bm I_{np_{t-1}}\}+\sum^{t}_{k=1}\mathscr{A}_{t-1}^{t-1-k}\text{vec}\{\widetilde\Phi_k(t-1)\widetilde\Phi^{\top}_k(t-1)\}\Big)\nonumber\\
				&~~~~~+\mathbb{I}_{\{p_t > p_{t-1}\}}\text{vec}\{P^{*-1}_{t-1}\} +\text{vec}\{\Phi_t\Phi^{\top}_t\}\Big\}.\label{421}
			\end{align}
			By the definition of $P^{*-1}_{t-1}$ and (\ref{vec2}), we have 
			\begin{align}
				&\text{vec}\{P^{*-1}_{t-1}\}\nonumber\\
				=&\text{col}\left\{\begin{pmatrix}
					P^{-1}_{t-1,1} & 0 \\
					0 & \bm I_{p_t-p_{t-1}}
				\end{pmatrix},\cdots,\begin{pmatrix}
					P^{-1}_{t-1,n} & 0 \\
					0 & \bm I_{p_t-p_{t-1}}
				\end{pmatrix}\right\}\nonumber\\
					=&\text{col}\Bigg\{
					\bm I_{p_{t}}+\sum^{t}_{k=1}\sum^n_{j=1}a^{(t-1-k)}_{1j}\widetilde\varphi_{k,j}(t)\widetilde\varphi^{\top}_{k,j}(t),\cdots \bm I_{p_{t}}+\sum^{t}_{k=1}\sum^n_{j=1}a^{(t-1-k)}_{nj}\widetilde\varphi_{k,j}(t)\widetilde\varphi^{\top}_{k,j}(t)
					\Bigg\}\nonumber\\
					=&\text{vec}\{\bm I_{np_{t}}\}+\sum^{t}_{k=1}\mathscr{A}_{t}^{t-1-k}\text{vec}\{\widetilde\Phi_k(t)\widetilde\Phi^{\top}_k(t)\}.\label{43}
				\end{align}
				Note that $\Phi_t=\Phi_t(t)$ and $\mathscr{A}_{t}$ is a stochastic matrix for all $t$, substituting (\ref{43}) into (\ref{421}),  it can be obtained that
				\begin{align}
					\text{vec}\{P^{-1}_{t}\}
					=&{\mathscr{A}_t}\Bigg\{\mathbb{I}_{\{p_t = p_{t-1}\}}\Big(\text{vec}\{\bm I_{np_{t-1}}\}
					+\sum^{t}_{k=1}\mathscr{A}_{t-1}^{t-1-k}\text{vec}\{\widetilde\Phi_k(t-1)\widetilde\Phi^{\top}_k(t-1)\}\Big)\nonumber\\
					&+\mathbb{I}_{\{p_t > p_{t-1}\}}\Big(\text{vec}\{\bm I_{np_{t}}\}+\sum^{t}_{k=1}\mathscr{A}_{t}^{t-1-k}\text{vec}\{\widetilde\Phi_k(t)\widetilde\Phi^{\top}_k(t)\}\Big) 
					+\text{vec}\{\Phi_t\Phi^{\top}_t\}\Bigg\}\nonumber\\
					=&\Bigg\{\mathbb{I}_{\{p_t = p_{t-1}\}}\Big({\mathscr{A}_t}\text{vec}\{\bm I_{np_{t-1}}\}+\sum^{t}_{k=1}\mathscr{A}_{t}^{t-k}\text{vec}\{\widetilde\Phi_k(t-1)\widetilde\Phi^{\top}_k(t-1)\}\Big)\nonumber\\
					&~~~\!+\!\mathbb{I}_{\{p_t > p_{t-1}\}}\Big({\mathscr{A}_t}\text{vec}\{\bm I_{np_{t}}\}\!+\!\sum^{t}_{k=1}\mathscr{A}_{t}^{t-k}\text{vec}\{\widetilde\Phi_k(t)\widetilde\Phi^{\top}_k(t)\}\Big) \Bigg\}\nonumber\\
					=&{\mathscr{A}_t}\text{vec}\{\bm I_{np_{t}}\}+\sum^{t}_{k=1}\mathscr{A}_{t}^{t-k}\text{vec}\{\widetilde\Phi_k(t)\widetilde\Phi^{\top}_k(t)\}\nonumber\\
					=&\text{vec}\{\bm I_{np_{t}}\}+\sum^{t}_{k=1}\mathscr{A}_{t}^{t-k}\text{vec}\{\widetilde\Phi_k(t)\widetilde\Phi^{\top}_k(t)\},
				\end{align}
				which concludes equation (\ref{vec}). 
				Thus for any $t\geq D_{\mathcal{G}}$, by  Remark \ref{RemaGrap},
				\begin{align}
					P^{-1}_{t,i}&=\bm I_{p_{t}}+\sum^{t}_{k=1}\sum^n_{j=1}a^{(t-k)}_{ij}\widetilde\varphi_{k,j}(t)\widetilde\varphi^{\top}_{k,j}(t)\label{s}\\
					&\geq \bm I_{p_{t}}+\sum^{t-D_{\mathcal{G}}}_{k=1}\sum^n_{j=1}a^{(t-k)}_{ij}\widetilde\varphi_{k,j}(t)\widetilde\varphi^{\top}_{k,j}(t)\geq a_{\min}\left\{\bm I_{p_{t}}+\sum^{t-D_{\mathcal{G}}}_{k=1}\sum^n_{j=1}\widetilde\varphi_{k,j}(t)\widetilde\varphi^{\top}_{k,j}(t)\right\}
				\end{align}	
				holds. From this, we have
				\begin{align}
					\lambda_{\min}\{P^{-1}_t\}\geq a_{\min}\lambda_{\min}\left\{\bm I_{p_{t}}+\sum^{t-D_{\mathcal{G}}}_{k=1}\sum^n_{j=1}\widetilde\varphi_{k,j}(t)\widetilde\varphi^{\top}_{k,j}(t)\right\}.\label{641}
				\end{align}	
				By the fact $a^{(t-k)}_{ij}\leq 1$ for all $i,j$, we have
				\begin{align}
					\lambda_{\max}\{P^{-1}_t\}=\max_{i\in\{1,...,n\}}	\lambda_{\max}\{P^{-1}_{t,i}\} \leq r_t.\label{471}
				\end{align}
				Then by Lemma \ref{lemma4.3} and (\ref{641})-(\ref{471}), it follows that
				\begin{align}
					&\!\!\operatorname{tr}(\widetilde{\Theta}^{\top}_{t}\widetilde{\Theta}_{t})\leq\frac{\operatorname{tr}\left(\widetilde{\Theta}^{\top}_{t}{P}_{t}^{-1}\widetilde{\Theta}_{t}\right)}{\lambda_{\min}(P^{-1}_{t})}=O\!\left(\frac{1}{\lambda_{\min}(t)}\left[p_t\log r_t\!+\sum^n_{i=1}\sum^{t}_{k=1}\|\epsilon_{k,i}\|^2\right]\right)  \text{a.s.}.
				\end{align}
				Moreover, we have
				$
				\sum^{\infty}_{q=p_t+1}\operatorname{tr}\left(\left(\Theta^{[q]}\right)^{\top}\Theta^{[q]}\right)
				= O(\gamma_t).
				$
				
				Thus combining this with (\ref{52}),
				by the definition of ${\hat\Theta}'_{t,i}$,  we have
				\begin{align*}
					\!\!&	\|{\hat\Theta}'_{t,i}-\Theta\|^2\leq \operatorname{tr}\left(({\hat\Theta}'_{t,i}-\Theta)^{\top}({\hat\Theta}'_{t,i}-\Theta)\right)\nonumber\\
					\!\!\leq& \operatorname{tr}\left(({\hat\Theta}_{t,i}-\Theta(t))^{\top}({\hat\Theta}_{t,i}-\Theta(t))\right)+\sum^{\infty}_{q=p_t+1}\operatorname{tr}\left(\left(\Theta^{[q]}\right)^{\top}\Theta^{[q]}\right)\nonumber\\
					\!\!=&O\!\left(\frac{1}{\lambda_{\min}(t)}\left[p_t\log r_t\!+\sum^n_{i=1}\sum^{t}_{k=1}\|\epsilon_{k,i}\|^2\right]+\gamma_t\right)  \text{a.s..}
				\end{align*}
				This completes the proof of the theorem.		
			\end{proof}	
           \begin{remark}
The upper bound in Theorem \ref{th1} reveals how excitation, dimension growth and truncation affect the estimation performance.
The denominator \(\lambda_{\min}(t)\) measures the overall excitation strength of the distributed system. It reflects the  advantage of distributed cooperation: by collecting spatio-temporal excitation information across nodes, the network can achieve sufficient excitation even when individual nodes cannot achieve it alone.
On the numerator, the term \(p_t \log r_t\) captures the estimation variance induced by dimension expansion: faster \(p_t\) growth adds more parameters to estimate from noisy observations, and zero-initialized newly activated dimensions introduce extra convergence error.
The term \(\gamma_t\) (or \(\epsilon_{k,i}\)) represents the model truncation bias, arising from truncating the true infinite-dimensional parameter to its first \(p_t\) components. Slower dimension growth leaves more high-order components truncated, resulting in a larger persistent bias.
Hence, this bound exhibits the classic bias-variance trade-off in selecting the dimension growth rate: accelerating dimension growth reduces truncation bias but increases estimation variance.
           \end{remark}

			Under some additional conditions, the above results in Theorem \ref{th1} can be  refined and articulated in a more concise form. Thus we have the following results.
			\begin{corollary}\label{re5.1}
				Under the conditions of Theorem \ref{th1}, if  the following conditions hold:
				
				(i) The unknown parameters satisfy $\Theta^{[q]}=O(\lambda^q)$ with some $0<\lambda<1$  for $q\geq 1$; 
				
				(ii) Each element of the regression vectors $\left\{\phi_{k,i}^{[q]}\right\}_{k\geq 1}$ is bounded for all $q \geq 1$ and $i\in\{1,...,n\}.$
				
				Then  by taking the regression lag $p_t=[\log^{\alpha} (t+e)]$ with $\alpha>1$ or $p_t=c[\log (t+e)]$ with $c\geq \frac{-1}{2\log \lambda}$ , we have
				\begin{align}
					\|{\hat\Theta}'_{t,i}-\Theta\|^2
					=O\left(\frac{p_t\log t}{\lambda_{\min}(t)}\right) \text{a.s.,}~\label{r3}
				\end{align}
				where ${\hat\Theta}'_{t,i}$ and $\lambda_{\min}(t)$ are defined in Theorem \ref{th1}.
			\end{corollary}
			\begin{proof}
				By the definition of $\epsilon^{\top}_{k,i}$ and the conditions (i) and (ii),  we have for $p_t=[\log^{\alpha} (t+e)]$ with $\alpha>1$ or $p_t=c[\log (t+e)]$ with $c> \frac{-1}{2\log \lambda}$ 
				\begin{align}
					\sum^n_{i=1}\sum^t_{k=1}\|\epsilon_{k,i}\|^2
					=&\sum^n_{i=1}\sum^t_{k=1}\left\|\sum^{\infty}_{q=p_{k-1}+1}\phi^{[q]}_{k,i}\Theta^{[q]}\right\|^2\nonumber\\
					\leq&\sum^n_{i=1}\sum^t_{k=1}\left(\sum^{\infty}_{q=p_{k-1}+1}\left\|\phi^{[q]}_{k,i}\right\|\left\|\Theta^{[q]}\right\|\right)^2
					=O\left(\sum^t_{k=1}\lambda^{2(p_{k-1})}\right)=O(1).
					\label{r4}
				\end{align}
				For $p_t=\frac{-1}{2\log \lambda}[\log (t+e)]$, it can be obtained from (\ref{r4}) that 
				\begin{align}
					\!\!\!\!\sum^n_{i=1}\sum^t_{k=1}\|\epsilon_{k,i}\|^2
					\!=\!O\left(\sum^t_{k=1}\lambda^{2(p_{k-1})}\right)\!=\!O(\log t)=O(p_t).
				\end{align}
				Moreover, by condition (ii) and the definition of $\lambda_{\min}(t)$, we have $\lambda_{\min}(t)=O(t)$.
				Combining this with (\ref{r4}), we can get 
				\begin{align}
					\gamma_t=O\left(\frac{\gamma_t t}{\lambda_{\min}(t)}\right)\label{70}
				\end{align}
				with $\gamma_t$ being defined in Theorem \ref{th1}.
				By condition (i) and  the selection of $p_t$, it can be obtained that
				$\gamma_t t=O(\lambda^{2(p_t+1)} t)=O(p_t\log t)$. By condition (ii) and the definition of $r_t$ in Theorem \ref{th1}, we have $r_t=O\left(\sum^t_{k=1}p^2_k\right)=O(p^2_tt)$.
				Thus by (\ref{r4})-(\ref{70}), we can see from Theorem \ref{th1} that 
				\begin{align*}
					\|{\hat\Theta}'_{t,i}-\Theta\|^2
					=O\left(\frac{p_t\log t}{\lambda_{\min}(t)}\right)~ \text{a.s..}
				\end{align*}
				This completes the proof of the theorem.
			\end{proof}
			
			\begin{remark}\label{re4.3}
				Suppose that Assumptions \ref{AssuSPR} and \ref{AssuGrap} are satisfied. If the unknown parameters satisfy $\Theta^{[q]}=O(q^{-a})$ with $a>1$, then following the proof of Corollary \ref{re5.1},   by taking the regression lag $p_t=O(t^{\alpha_1})$ with $\alpha_1\geq \frac{1}{2a-1}$, we also have
				\begin{align*}
					\|{\hat\Theta}'_{t,i}-\Theta\|^2
					=O\left(\frac{p_t\log t}{\lambda_{\min}(t)}\right) \text{a.s..}
				\end{align*}
			\end{remark}
			The proof is similar to Corollary \ref{re5.1} and we omit it here.

			\begin{remark}\label{re5.2}
				From (\ref{r3}), the strong convergence of Algorithm \ref{algorithm} is guaranteed provided the following cooperative excitation condition is satisfied:
				\begin{align}
					\frac{p_t\log t}{\lambda_{\min}(t)}\xrightarrow{t\to\infty}0,~ \text{a.s..}\label{e17}
				\end{align}  
				Some remarks on this condition are provided below.
				
				(i) \textbf{Finite-dimensional extension:} When the system model (\ref{model}) reduces to a finite-dimensional one, i.e., $p_t$ is taken as a  constant $p$, condition (\ref{e17})  simplifies to the weakest possible cooperative non-persistent excitation (PE) condition for finite-dimensional regression models given in \cite{Xie2021}, namely
				\begin{align}
					\frac{\log t}{\lambda_{\min}\left(\sum^n_{i=1}\sum^{t-D_{\mathcal{G}}}_{k=1}\varphi_{k,i}\varphi^{\top}_{k,i}\right)}\xrightarrow{t\to\infty}0,~ \text{a.s..}\label{e18}
				\end{align} 
				Thus, the result in \cite{Xie2021} appears as a special case of (\ref{r3}) when $p_t$ is a constant. This shows that condition (\ref{e17}) constitutes a natural extension of condition (\ref{e18}) from the finite-dimensional to the infinite-dimensional setting.
				
				(ii) \textbf{Infinite-dimensional comparison:} 
				For choices such as $p_t=c[\log (t+e)]$ with sufficiently large $c$  or  $p_t=[\log^{\alpha} (t+e)]$ with $\alpha>1$, or even  faster-growing  $p_t= [t^{\alpha_1}]$ with $0 < \alpha_1 < 1$, the condition (\ref{e17}) is much weaker than the standard cooperative PE condition, i.e.,  $t=O(\lambda_{\min}(t)),~ \text{a.s..}$
				Condition (\ref{e17}) also generalizes the non-cooperative excitation condition from \cite{Guo1991,DAI2025, DAI20251} (i.e., $\mathcal{A} = I_{n}$) by relaxing the requirement on individual agents $i$:
				\begin{gather}\label{noncoo}
					\frac{p_t\log t}{\lambda_{\min}\left(\sum^{t}_{k=1}\widetilde\varphi_{k,i}(t)\widetilde\varphi^{\top}_{k,i}(t)+ I_{p_t}\right)}\xrightarrow{t\to\infty}0.
				\end{gather}
				This implies that even when no single agent can accurately estimate the unknown parameter matrix using traditional non-cooperative algorithms, the multi-agent system as a whole may still accomplish the estimation task through the distributed algorithm, as numerically validated in Section \ref{simulation}.
			\end{remark}

			\subsection{Regret analysis}
			We now proceed to analyze the accuracy of the online prediction scheme. Since equation (\ref{y1}) holds for all $k > 0$, we take conditional expectations on both sides under Assumption \ref{AssuSPR}, which yields the best  prediction for ${y}_{k+1,i}^{\top}$ of agent $i$ in the mean square sense:
			\begin{gather}
				\mathbb{E}({y}_{k+1,i}^{\top} \mid \mathscr{F}_k) = \mathbb{E}(\bar\varphi^{\top}_{k,i} \Theta(k-1) + \epsilon^{\top}_{k,i} \mid \mathscr{F}_k). \label{53}
			\end{gather}
			By substituting the unknown parameter $\Theta(k-1)$ in (\ref{53}) with agent $i$'s estimate $\hat{\Theta}_{k-1,i}$ and neglecting the residual term $\epsilon^{\top}_{k,i}$, we obtain the following online predictor for ${y}_{k+1,i}^{\top}$:
			\begin{gather}
				\hat{y}_{k+1,i}^{\top} = \mathbb{E}(\bar\varphi^{\top}_{k,i} \hat{\Theta}_{k-1,i} \mid \mathscr{F}_k). \label{55}
			\end{gather}
			The difference between the ideal prediction and the online predictor is commonly termed the ``regret'' \cite{Xie2021}. Accordingly, we define the following regret for agent $i$ at time $k$ as:
			\begin{align*}
				R_{k,i} = \left\| \mathbb{E}({y}_{k+1,i}^{\top} \mid \mathscr{F}_k) - \hat{y}_{k+1,i}^{\top} \right\|^2.
			\end{align*}

			Although the influence of noise makes it impractical to expect the per-step regret to be consistently small, we establish that the  regret is small in an average sense. Theorem \ref{th3} demonstrates that the  accumulated regret from $k=1$ to $t$ is of order $o(t)$.
			\begin{theorem}\label{th3}
				Consider the model (\ref{model})-(\ref{con}) and Algorithm \ref{algorithm}.  Under Assumptions \ref{AssuSPR} and \ref{AssuGrap}, if  $\left\{\phi_{k,i}^{[q]}\right\}_{k\geq 1}$ is bounded for every $q \geq 1$ and $i\in\{1,...,n\}$, then the accumulated regret over networks has the following upper bound as $t\rightarrow\infty$:
				\begin{align*}
					\sum^n_{i=1}\!\sum^{t}_{k=1}\!R_{k,i}
					\!=O\left(p^2_t\log t+p_t\sum^n_{i=1}\sum^{t}_{k=1}\|\epsilon_{k,i}\|^2\right) ~{\rm a.s..}
				\end{align*}
			\end{theorem}
			\begin{proof}
				By Assumption \ref{AssuSPR}, it is clear that $\bar\varphi_{k,i}, \epsilon^{\top}_{k,i}$ and $\hat\Theta_{k-1,i}$ is $\mathscr{F}_k$-measurable. Thus from (\ref{53}) and (\ref{55}), we can obtain that
				\begin{align}
					&\sum^n_{i=1}\sum^{t}_{k=1}R_{k,i}=\sum^n_{i=1}\sum^{t}_{k=1}\|\mathbb{E}({y}_{k+1,i}^{\top}|\mathscr{F}_k)-\hat{y}_{k+1,i}^{\top}\|^2\nonumber\\
					=&\sum^n_{i=1}\sum^{t}_{k=1}\|\bar\varphi^{\top}_{k,i}{\widetilde\Theta_{k-1,i}}+\epsilon^{\top}_{k,i}\|^2
					\leq 2\sum^n_{i=1}\sum^{t}_{k=1}\|\bar\varphi^{\top}_{k,i}{\widetilde\Theta_{k-1,i}}\|^2+2\sum^n_{i=1}\sum^{t}_{k=1}\|\epsilon_{k,i}\|^2.\label{577}
				\end{align}
				Note that   $\lambda_{\max}\big(P_{k+1,i}\big) \!=\!\frac{1}{\lambda_{\min}\big(P^{-1}_{k+1,i}\big)}\!\leq\!1$ 
				holds for all $k$ and $i$ by (\ref{vec2}). Then by (\ref{14}) and (\ref{31}) and the boundedness of $\left\{\phi_{k,i}^{[q]}\right\}_{k\geq 1}$,  we have
				\begin{align*}
					b^{-1}_k =\bm I_{n} \!+\! \bar\Phi_{k}^{\top}P_{k-1}\bar\Phi_{k}+\!  \mathbb{I}_{\{p_k > p_{k-1}\}}{\bar{\bar\Phi}}_{k}^{\top}{\bar{\bar\Phi}}_{k} =O(p_k),
				\end{align*}	
				which implies that
				\begin{align}
					b^{-1}_{k,i}=O(p_k).\label{57}
				\end{align}
				Combining (\ref{57}) with Lemma \ref{lemma4.3}, we have
				\begin{align}
					\sum^n_{i=1}\sum^{t}_{k=1}\|\bar\varphi^{\top}_{k,i}{\widetilde\Theta_{k-1,i}}\|^2
					=&\sum^n_{i=1}\sum^{t}_{k=1}b^{-1}_{k,i}b_{k,i}\|\bar\varphi^{\top}_{k,i}\widetilde{\Theta}_{k,i}\|^2\nonumber\\
					=&O\left(p^2_t\log^+\{ \lambda_{\max}(P^{-1}_t)\}+p^2_t+p_t\sum^n_{i=1}\sum^{t}_{k=1}\|\epsilon_{k,i}\|^2\right)\nonumber\\
					=&O\left(p^2_t\log t+p_t\sum^n_{i=1}\sum^{t}_{k=1}\|\epsilon_{k,i}\|^2\right) ~{\rm a.s.,}\label{79}
				\end{align}
				where equation (\ref{s}) and the  boundedness of $\{\phi_{k,i}^{[q]}\}_{k\geq 1}$ are used in the last inequality.
				Combining (\ref{577}) and (\ref{79}), we complete the proof of the theorem.
			\end{proof}
			
			\begin{remark}\label{re4.5}
				Under the  conditions of Theorem \ref{th3} or Remark \ref{re4.3}, we can obtain the  following bound of the accumulated regret over networks in a more concise form:
				\begin{align}
					\sum^n_{i=1}\!\sum^{t}_{k=1}\!R_{k,i}
					\!=O\left(p^2_t\log t\right) ~{\rm a.s..}\label{80}
				\end{align}
				Thus for choices such as $p_t=c[\log (t+e)]$ with sufficiently large $c$ or $p_t=[\log^{\alpha} (t+e)]$ with $\alpha>1$, or   $p_t= [t^{\alpha_1}]$ with $0 < \alpha_1 < \frac{1}{2}$, we have $\frac{\sum^n_{i=1}\!\sum^{t}_{k=1}\!R_{k,i}}{nt}\to 0$ as $t\to\infty$~{\rm a.s..}
				Furthermore, when $p_t$ is taken as the upper bound of the finite-order system dimension, (\ref{80}) coincides with the classical convergence properties of the distributed least squares algorithm for finite-dimensional stochastic regression models. In this case, the asymptotic order of the accumulative regret, $O(\log t)$, may be optimal for adaptive prediction \cite{Xie2021}.
			\end{remark}
			

			Compared with the theoretical results in \cite{Gan2025CCC}, the algorithm proposed in this paper optimizes the convergence rate of the prediction ``regret" without imposing any other  conditions on estimation convergence. This improvement is primarily attributed to the adaptive mechanism of Algorithm \ref{algorithm}, which eliminates the need to consider the consistency issue between the current time step $t$ and the parameter estimation iteration count $k$. In contrast, the proposed algorithm can be more conveniently applied to adaptive feedback control systems while effectively reducing computational complexity.

				\section{Simulation Results}\label{simulation}
				This section presents numerical simulations to verify the theoretical convergence, cooperative estimation mechanism, and effectiveness of the adaptive dimension-growth strategy in the proposed Algorithm \ref{algorithm}. 
				
				We consider a network consisting of \(n=18\) agents. Each agent performs an online estimation task based on the linear regression model (\ref{model}), where \(m=2\) denotes the output dimension, and 
				$$\Theta^\top = \begin{pmatrix}
					0.5 & 0.5^2 & 0.5^3 & \dots 0.5^{100}\\
					0.6 & 0.6^2 & 0.6^3 & \dots 0.6^{100}
				\end{pmatrix}\in \mathbb{R}^{2\times 100} $$
				represents the unknown parameter matrix truncated to a finite maximum dimension $i.e., 100$. The colored non-independent and identically distributed (non-i.i.d.) regression vector \(\phi_{k,i}\) is constructed as follows: for agent \(i\) at time \(k\), its \(q\)-th element (\(q\neq i\)) is defined as \(\phi_{k,i}^{[q]} = \sum_{j=1}^k o_{j,i}\), where \(\{o_{j,i}\}\) follows the uniform distribution \(\mathcal{U}(-i,i)\) over the interval \([-i, i]\); the \(i\)-th element is fixed to zero to remove self-excitation. This setting ensures that no single agent can independently satisfy the excitation condition (\ref{noncoo}), while the entire multi-agent network collectively fulfills the cooperative excitation condition (\ref{e17}) to achieve consistent estimation. The measurement noise \(w_{k,i}\) is an independent zero-mean Gaussian noise following \(\mathcal{N}(0, 0.4^2, m, 1)\), where the notation  \(\mathcal{N}(\mu, \sigma, s, r)\)
				denotes an $s\times r$ matrix with each element obeying the Gaussian
				distribution with mean $\mu$ and standard deviation $\sigma$.
				
				The communication topology is depicted in Fig. \ref{graph}, which is an undirected connected graph. Its doubly stochastic weight matrix is designed using the Metropolis rule in \cite{new05}. The adaptive dimension-growth sequence in Algorithm \ref{algorithm} is set to \(p_t = \lfloor \log^{2}(t+e) \rfloor\) with an initial dimension \(p_0=1\), which satisfies the slow-growth requirement in the theoretical analysis. We perform \(S=50\) independent experiments to reduce random fluctuations.
				
				1) We compare our algorithm with the non-cooperative  LS algorithm, where the weighted adjacency matrix is set to the identity matrix \(\bm I_n\) (i.e., no communication between agents). For initialization, the local parameter estimates are randomly set as \(\hat{\Theta}_{0,i} \sim 5\mathcal{N}(0, 1,p_0,m)\), and the information matrix is initialized to the positive definite matrix \(P_{0,i} = \bm I_{p_0}\). Accordingly, we obtain \(S=50\) sequences for each agent:
				$$\left\{\|\hat{\Theta}_{t,i}^{(s)} - \Theta_{1:p_t}\|\right\}_{t\geq 0},~ s=1,...,50,~ i=1,...,18,$$
				where \(\hat{\Theta}_{t,i}^{(s)}\) denotes the estimate of agent \(i\) in the \(s\)-th experiment, and \(\Theta_{1:p_t}\) is the true parameter truncated to dimension \(p_t\). We define 
				\begin{gather}
					e_{t,i} = \frac{1}{S}\sum^S_{s=1}\|\hat{\Theta}_{t,i}^{(s)} - \Theta_{1:p_t}\|\label{ess}
				\end{gather}
				as the average estimation error of agent \(i\), as plotted in Fig. \ref{fig:coo_vs_noncoo}. This figure shows the estimation errors of all agents under non-cooperative and cooperative frameworks. The non-cooperative LS algorithm produces constant steady-state errors because no single agent can satisfy the excitation condition (\ref{noncoo}) mentioned earlier. In contrast, when all agents cooperate with their neighbors using Algorithm \ref{algorithm} over the connected network, the estimation error converges to zero faster and remains stable. This result verifies the cooperative advantage of multi-agent systems, as described in Remark \ref{re5.2}.
				
				2) We also compare our algorithm with the fixed-dimension distributed least squares (Fixed-DLS) algorithm with fixed dimensions \(p=2,3,4,5\), and the results are presented in Fig. \ref{fig:adaptive_vs_fixed}. All algorithms adopt the same recursive update rule and only differ in dimension adjustment methods. The  initialization  of Fixed-DLS are set as \(\hat{\Theta}_{0,i} \sim 5\mathcal{N}(0, 1,p,m)\), and  \(P_{0,i} =  \bm I_{p}\) with $p=2,3,4,5$. The network-averaged estimation error is calculated as \(\bar{e}(t) = \frac{1}{n}\sum_{i=1}^n e_{t,i}\) with $e_{t,i}$ being defined in (\ref{ess}). As shown in Fig. \ref{fig:adaptive_vs_fixed},  for Fixed-DLS schemes, steady-state estimation error reduces as dimension $p$ increases, yet non-negligible steady-state bias persists. In contrast, the proposed algorithm achieves lower steady-state error, demonstrating its superiority over conventional fixed-dimension DLS methods in ultimate estimation accuracy. This confirms the effectiveness of the dimension-growth rule in matching the true infinite-dimensional parameter structure.


				\begin{figure}[htb!]
					\centering
					\includegraphics[width=0.5\textwidth]{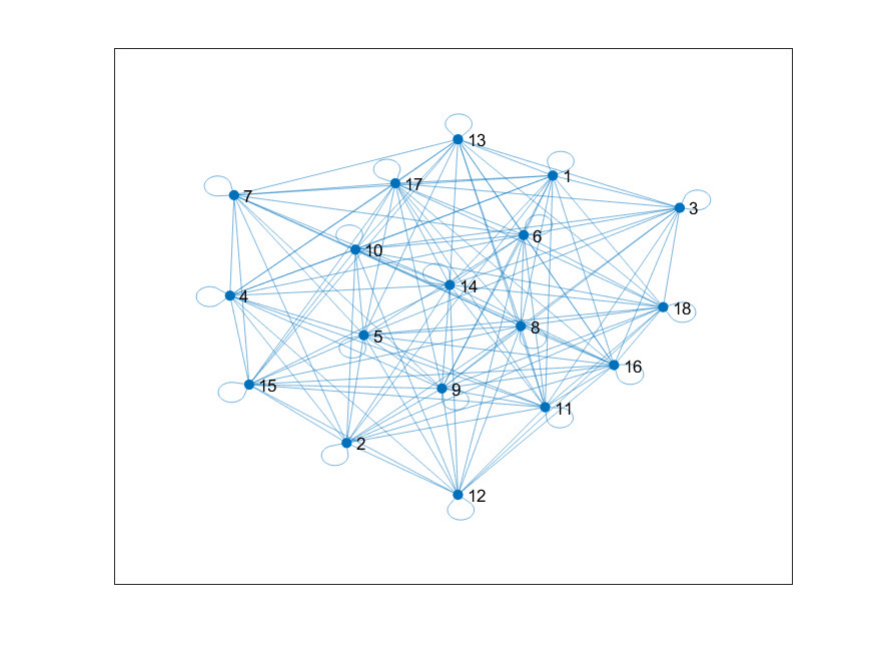}
					\caption{Communication topology among 18 agents}
					\label{graph}
				\end{figure}
				
				\begin{figure*}[htb!]
					\centering
					\includegraphics[width=1\textwidth]{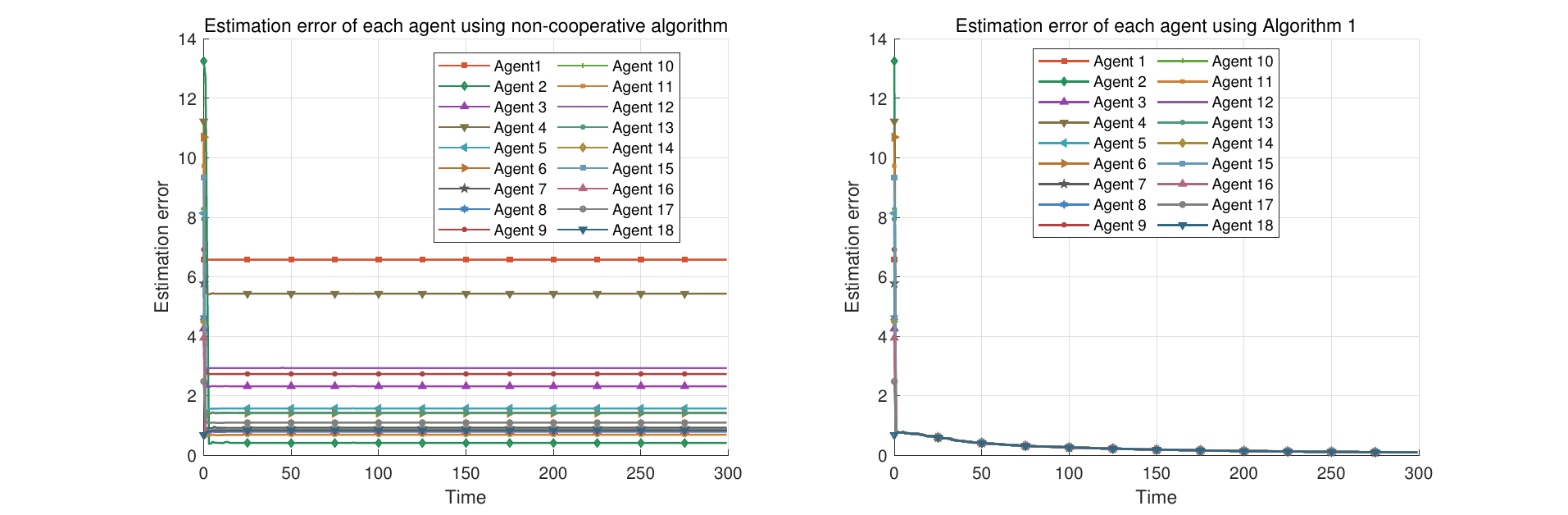}
					\caption{Estimation errors of all 18 agents: (left) non-cooperative LS algorithm; (right) proposed cooperative distributed algorithm (Algorithm \ref{algorithm}).}
					\label{fig:coo_vs_noncoo}
				\end{figure*}
				
				\begin{figure}[htb!]
					\centering
					\includegraphics[width=0.5\textwidth]{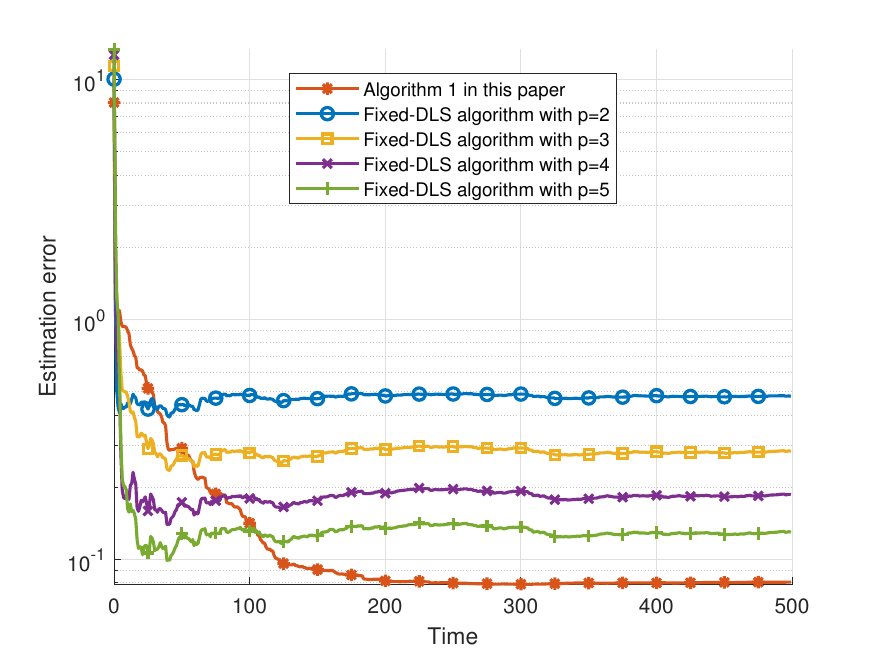}
					\caption{Network-averaged estimation errors of proposed Algorithm \ref{algorithm} and fixed-dimension distributed LS algorithms with \(p=2,3,4,5\).}
					\label{fig:adaptive_vs_fixed}
				\end{figure}

				\section{Concluding Remarks}\label{section_conclusion}
				This paper proposed a fully distributed online RLS algorithm with dynamically growing dimensions to estimate the infinite unknown parameters of a stochastic large regression model.
 The algorithm adopts a dimension-adaptive local RLS update scheme with zero-initialization augmentation for newly activated dimensions, and achieves consensus estimation  via a covariance-intersection-based neighbor information fusion rule.
				We analyzed the convergence properties of our algorithm under a cooperative non-persistent condition tailored for the random regression vectors, without requiring assumptions of independence, stationarity, or ergodicity. This relaxation enables our theoretical results to be applicable to feedback systems.
				Moreover, we established an almost sure upper bound of the accumulated prediction regret without any excitation conditions on the regressors. Our results not only recover the optimal finite-dimensional  estimation bound as a special case when the dimension is fixed, but also reveal the  trade-off between dimension expansion speed and estimation performance.
				 Future work included privacy protection and complexity reduction in large- model distributed learning, the fusion of distributed adaptive estimation and control, estimation of nonlinear systems, and extensions to non-ideal communication scenarios.
				
					
					

				\bibliographystyle{siamplain}
				\bibliography{my_bib}
			\end{document}